\documentclass[11pt]{article}

\usepackage{amsmath,amssymb,amsthm,mathrsfs}
\usepackage{enumitem}
\usepackage{geometry}
\usepackage{xcolor}
\usepackage{url}
\numberwithin{equation}{section}

\newtheorem{theorem}{Theorem}[section]
\newtheorem{lemma}[theorem]{Lemma}
\newtheorem{proposition}[theorem]{Proposition}
\newtheorem{corollary}[theorem]{Corollary}
\newtheorem{remark}[theorem]{Remark}

\newcommand{\PP}{\mathbb{P}}
\newcommand{\BE}{\mathbb{E}}
\newcommand{\E}{\mathrm{e}}
\newcommand{\D}{\mathrm{d}}

\newcommand{\vertk}{\stackrel{{\cal D}}{\longrightarrow}}
\newcommand{\stk}{\stackrel{{\mathbb{P}}}{\longrightarrow}}
\newcommand{\edist}{\stackrel{{\cal D}}{=}}
\newcommand{\argmin}{\operatorname*{arg\,min}}
\newcommand{\argmax}{\operatorname*{arg\,max}}

\title{Joint extreme-value asymptotics for the radius and centre\\
of smallest enclosing balls}

\author{Norbert Henze\footnote{Institute of Stochastics, Karlsruhe Institute of Technology (KIT), 76131 Karlsruhe, Germany. e-mail:
henze@kit.edu}} %\ and Ilya Molchanov\footnote{Institut f\"ur Mathematische Statistik und %Versicherungslehre. Alpeneggstrasse 22, 3012 Bern, Switzerland. e-mail: ilya.molchanov@unibe.ch}}

\date{}

\begin{document}
\maketitle

\begin{abstract}
We study the joint extreme-value behaviour of the radius and centre of
the smallest ball containing an i.i.d.\ spherically symmetric sample.
In the Gumbel and Weibull regimes, we obtain joint limit theorems and
show that the enclosing ball is asymptotically supported by exactly
$d+1$ sample points, thereby resolving an open support problem in the
Gumbel case. In the Fr\'echet regime, the limit is given by the
smallest enclosing ball of an extremal Poisson cloud. For planar
Gaussian samples, we derive an explicit joint limit law from Daniels'
finite-sample formula and show that this law is universal throughout
the planar Gumbel domain. In the bounded-support Weibull regime, we additionally
obtain uniform tail bounds, moment convergence and an explicit
limiting density, while the uniform-ball case admits an interpretation
in terms of a Poisson zero cell.
\end{abstract}

\bigskip
\noindent\textbf{MSC 2020:} Primary 60G70; Secondary 60D05; 60F05.

\noindent\textbf{Keywords:} smallest enclosing ball; enclosing radius;
centre; extreme values; Poisson point process; stochastic geometry;
Gaussian sample; spherical symmetry.

%%%%%%%%%%%%%%%%%%%%%%%%%%%%%%%%%%%%%%%%%%%%%%%%%%%%%%%%%%%%%%%%%%%%%%%%%%%%%%
\section{Introduction}\label{sec:intro}
%%%%%%%%%%%%%%%%%%%%%%%%%%%%%%%%%%%%%%%%%%%%%%%%%%%%%%%%%%%%%%%%%%%%%%%%%%%%%%%

Let $X,X_1,X_2,\ldots$ be independent and identically distributed
random vectors in $\mathbb R^d$, $d\geq2$, with an absolutely continuous
spherically symmetric distribution.  Let
\[
 R_n:=\inf_{c\in\mathbb R^d}\max_{1\leq j\leq n}\|X_j-c\|
\]
denote the radius of the smallest enclosing ball of
$X_1,\ldots,X_n$, and let $C_n$ be its centre.  Thus
\[
 R_n=\max_{1\leq j\leq n}\|X_j-C_n\|.
\]
Throughout, $\|\cdot\|$ denotes the Euclidean norm on $\mathbb R^d$,
and $x\cdot y$ denotes the usual Euclidean inner product.
The centre is uniquely determined for every finite configuration:
if two distinct centres gave the same minimal radius, the
parallelogram identity would give a strictly smaller radius at their
midpoint.

Put $\rho=\|X\|$. Since the distribution of $X$ is absolutely
continuous, $\PP(X=0)=0$. On $\{X\neq0\}$ let
$U=X/\|X\|$, with an arbitrary definition on the null event
$\{X=0\}$. By spherical symmetry, $U$ is uniformly distributed on
\[
\mathbb S^{d-1}:=\{x\in\mathbb R^d:\|x\|=1\},
\]
is independent of $\rho$, and $X=\rho U$ almost surely. We use the
analogous notation $X_j=\rho_jU_j$ for the sample. We denote by
$\sigma_{d-1}$ the normalized surface measure on $\mathbb S^{d-1}$,
so that $\sigma_{d-1}(\mathbb S^{d-1})=1$. Moreover, let
\[
\mathbb B^d:=\{x\in\mathbb R^d:\|x\|\le1\}
\]
denote the closed unit ball, and write $\kappa_d$ for its $d$-dimensional
Lebesgue measure.

Henze and Klatt~\cite{hekl} studied the asymptotic behaviour of $R_n$
in the three classical extreme-value regimes. Their results describe
the limiting radius through Poisson point processes and geometric
variational problems, and also provide tightness results for the
corresponding centres. The purpose of the present paper is to
strengthen this theory by deriving joint limit theorems for radius and
centre. In contrast to Henze and Klatt~\cite{hekl}, we identify the
limiting centre, prove joint convergence of radius and centre, and
establish the asymptotic supporting-simplex structure.

The sample diameter
\[
D_n:=\max_{1\le i<j\le n}\|X_i-X_j\|
\]
has been studied quite extensively in related settings. For
spherically symmetric distributions, see Jammalamadaka and
Janson~\cite{jamsva15}; for bounded-support models, see
Lao~\cite{lao10}. Further results for elliptical distributions in the
Gumbel domain were obtained by Demichel, Fermin and Soulier~\cite{dfs15}.
Compared with $D_n$, considerably less has been known about the joint
behaviour of the smallest enclosing radius $R_n$ and its centre $C_n$.

In the Gumbel and Weibull regimes, the limiting radius is the optimal
value of a random Poisson variational problem.  We show that the
corresponding optimizer is almost surely unique and obtain joint
convergence of the normalized radius and centre to the value--optimizer
pair.  The limiting optimizer is supported by exactly $d+1$ active
Poisson constraints, which also yields asymptotic $(d+1)$-point support
of the sample enclosing ball.  In particular, this answers affirmatively the higher-dimensional Gumbel support problem posed in Henze and Klatt~\cite[Problem~7.6]{hekl}. In the Weibull case we obtain, in
addition, uniform tail bounds, moment convergence and an explicit joint
density with a gamma factorisation.

In the Fr\'echet regime, the limiting object is instead the smallest
enclosing ball of the extremal Poisson cloud.  We prove joint
convergence of its radius and centre by exploiting continuity of the
smallest-enclosing-ball map under Hausdorff convergence.

A classical result of Daniels~\cite{dan52} provides an explicit
benchmark.  For planar Gaussian samples, Daniels obtained an exact
finite-sample formula for the radius and the distance of the centre
from the origin.  We use this formula to derive their joint limiting
distribution.  Combined with the general Gumbel theorem, this gives an
explicit universal joint limit throughout the planar Gumbel
max-domain of attraction.

For samples uniformly distributed on the unit ball, the Weibull limit
is closely related to the random-set limit of Kabluchko, Marynych and
Molchanov~\cite{KMM}.  In this case the limiting radius--centre pair
can be identified with the inradius and the reflected incenter of a
Poisson zero cell.

The paper is organized as follows.  Section~\ref{sec:daniels} treats
the planar Gaussian case.  Section~\ref{sec:gumbel} establishes the
general Gumbel joint limit and the asymptotic supporting simplex.
Section~\ref{sec:frechet} deals with the Fr\'echet regime, and
Section~\ref{sec:weibull} treats the Weibull regime.

Throughout, $Y\edist Z$ denotes equality in distribution, whereas
$Y_n\vertk Y$ and $Y_n\stk Y$ denote convergence in distribution and
in probability, respectively.

%%%%%%%%%%%%%%%%%%%%%%%%%%%%%%%%%%%%%%%%%%%%%%%%%%%%%%%%%%%%%%%%%%%%%%%%%%%%%%%
%
% Section 2
\section{Daniels' exact formula and the joint Gaussian limit}\label{sec:daniels}
%
%%%%%%%%%%%%%%%%%%%%%%%%%%%%%%%%%%%%%%%%%%%%%%%%%%%%%%%%%%%%%%%%%%%%%%%%%%%%%%%

In this section, we assume that $X$ has a standard normal distribution in
$\mathbb R^2$. This special case was considered by Daniels~\cite{dan52}, who obtained an
exact joint density for $R_n$ and $\|C_n\|$.  We first recall the form
of his result that is most convenient for asymptotic analysis.

For $c\in\mathbb R^2$ and $r\geq 0$, let $B(c,r):=\{x\in\mathbb R^2:\|x-c\|\leq r\}$ 
denote the closed Euclidean ball with centre $c$ and radius $r$. 
For $r,p\geq 0$, let
\[
P(r,p):=\PP\{X_1\in B(c,r)\}, \qquad \|c\|=p.
\]
By rotational symmetry, this probability depends on $c$ only through
$\|c\|=p$, so that $P(r,p)$ is well defined.

Here and below, $I_\nu$ denotes the modified Bessel function of the first kind of order $\nu$;
see, e.g., \cite{DLMF}. In particular,
\[
I_0(x)=\frac{1}{2\pi}\int_0^{2\pi} \E^{x\cos\theta}\,\D \theta,
\qquad x\geq0,
\]
and $I_1=I_0'$.  Daniels \cite{dan52} showed that
\begin{equation}\label{eq:Daniels-expansion}
\begin{split}
 1-P(r,p)
 ={}&\exp\left\{-\frac{r^2+p^2}{2}\right\}
 \bigg\{ I_0(rp)+\frac{p}{r}I_1(rp)
       +\frac{p^2}{r^2}I_2(rp)+\cdots\bigg\}.
\end{split}
\end{equation}
The series converges for all $r,p>0$, since
$\left(\frac{p}{r}\right)^k I_k(rp)
\le I_0(rp)(p^2/2)^k/k!$.

Daniels separates the possible supporting configurations of the
smallest enclosing circle into two cases: either the circle is
supported by three points forming an acute triangle, or it is supported
by two diametrically opposite points.  For unit component variance,
his formulas (3.4)--(3.6) yield the following joint density of
$(R_n,\|C_n\|)$:
\begin{equation}\label{eq:Daniels-joint-density}
\begin{split}
 f_n(r,p)
={}&n(n-1)(n-2)
 \exp\left\{-\frac32(r^2+p^2)\right\}
 r^2 I_1(rp)P(r,p)^{n-3}\\
&\quad +2n(n-1)
 \exp\left\{-(r^2+p^2)\right\}
 rp\,P(r,p)^{n-2},
 \qquad r,p>0.
\end{split}
\end{equation}
The first term corresponds to three supporting points and the second to
two diametrically opposite supporting points.
We now extract the joint extreme-value limit from \eqref{eq:Daniels-joint-density}.  To this end, put
\[
%\begin{equation}\label{eq:gaussian-an-bn}
 a_n:=\sqrt{2\log n},\qquad b_n:=\frac1{a_n},
%\end{equation}
\]
and define
\[
%\begin{equation}\label{eq:Tn-Sn}
 T_n:=\frac{R_n-a_n}{b_n}       =a_n(R_n-a_n),  \qquad  S_n:=\frac{\|C_n\|}{b_n}   =a_n\|C_n\|.
%\end{equation}
\]

\smallskip
With this normalization, Daniels' formula yields the following joint
limit theorem.

\medskip

%%%%%%%%%%%%%%%%%%%%%%%%%%%%%%%%%%%%%%%%%%%%%%%%%%%%%%%%%%%%%%%%%%%%%%%%%%%%%%%%%%%%%%%%%%%%%%
%
% Theorem 2.1
%
\begin{theorem}[Joint Gaussian limit]\label{thm:Daniels-joint-limit}
%%%%%%%%%%%%%%%%%%%%%%%%%%%%%%%%%%%%%%%%%%%%%%%%%%%%%%%%%%%%%%%%%%%%%%%%%%%%%%%%%%%%%%%%%%%%%%%%%%%
As $n\to\infty$,
\[
 (T_n,S_n)\vertk(T,S),
\]
where $(T,S)$ has density
\begin{equation}\label{eq:joint-limit-density}
 f(t,s)  =\E^{-3t}I_1(s)   \exp\{-\E^{-t}I_0(s)\},  \qquad t\in\mathbb R,\ s>0.
\end{equation}
In fact, the densities of $(T_n,S_n)$ converge to $f$ in $L^1$.
\end{theorem}
%%%%%%%%%%%%%%%%%%%%%%%%%%%%%%%%%%%%%%%%%%%%%%%%%%%%%%%%%%%%%%%%%%%%%%%%%%%%%%%%%%%%%%%%%%%%%%%%%%%

\smallskip

\begin{proof}
For fixed $t\in\mathbb R$ and $s>0$, put $ r_n:=a_n+b_nt$, $ p_n:=b_ns$. Then
\begin{equation}\label{eq:basic-gaussian-expansions}
 r_n^2=a_n^2+2t+o(1),  \qquad  r_np_n\rightarrow s,  \qquad  \frac{p_n}{r_n}\rightarrow0.
\end{equation}
Since $a_n^2=2\log n$, we have
\begin{equation}\label{eq:gaussian-tail-scale}
 n\exp\left\{-\frac{r_n^2}{2}\right\}  \rightarrow \E^{-t}.
\end{equation}
%%%%%%%%%%%%%%%%%%%%%%%%%%%%%%%%%%%%%%%%%%%%%%%%%%%%%%%%%%%%%%%%%%%%%%%%%%%%%%%%%%%%%%%%%%%%%%%%%
Put $q_n:=p_n/r_n$. Since $r_np_n\to s$ and $q_n\to0$, the terms
following $I_0(r_np_n)$ in \eqref{eq:Daniels-expansion} are $o(1)$.
Indeed, for some $K>s$ and all sufficiently large $n$, using
$I_k(x)\le I_0(x)$ for $x\ge0$, we have 
\[
 \sum_{k\ge1} q_n^k I_k(r_np_n)
 \le I_0(K)\frac{q_n}{1-q_n},
\]
and the right-hand side converges to $0$. Together with \eqref{eq:basic-gaussian-expansions}
and \eqref{eq:gaussian-tail-scale}, this yields
\begin{equation}\label{eq:P-asymptotic}
 n\{1-P(r_n,p_n)\}
 \rightarrow \E^{-t}I_0(s).
\end{equation}
Consequently, for each fixed integer $j$,
\begin{equation}\label{eq:P-power-asymptotic}
 P(r_n,p_n)^{n-j}  \rightarrow  \exp\{-\E^{-t}I_0(s)\}.
\end{equation}
Let $f_n$ denote the joint density of $(R_n,\|C_n\|)$ given in \eqref{eq:Daniels-joint-density}, and let $g_n$ denote the joint density of $(T_n,S_n)$. By the change of variables
$r=a_n+b_n t$ and $p=b_n s$, we have 
$g_n(t,s)=b_n^2 f_n(r_n,p_n)$.

Consider first the contribution from the first term in
\eqref{eq:Daniels-joint-density}.  In view of
\eqref{eq:basic-gaussian-expansions}--\eqref{eq:P-power-asymptotic},
\begin{align*}
&b_n^2 n(n-1)(n-2)
 \exp\left\{-\frac32(r_n^2+p_n^2)\right\}  r_n^2I_1(r_np_n)P(r_n,p_n)^{n-3}\\
&\hspace{35mm}\rightarrow  \E^{-3t}I_1(s)  \exp\{-\E^{-t}I_0(s)\}.
\end{align*}
For the second term, we have $2b_n^2 n(n-1) \exp\{-(r_n^2+p_n^2)\} r_np_nP(r_n,p_n)^{n-2} \rightarrow0$,
since $n^2\exp(-r_n^2)\to\E^{-2t}$ while
$b_n^2r_np_n\to0$.  Hence $g_n(t,s)\rightarrow f(t,s)$ for every $t\in\mathbb R$ and $s>0$.

It remains to check that $f$ integrates to one.  Put $w=\E^{-t}$.
Since $I_0'(s)=I_1(s)$ and $I_0(0)=1$,
\[
 \int_0^\infty I_1(s)\E^{-wI_0(s)}\,\D s
 =\int_1^\infty \E^{-wy}\,\D y
 =\frac{\E^{-w}}{w}.
\]
Therefore
\[
 \int_{-\infty}^{\infty}\int_0^\infty f(t,s)\,\D s\,\D t
 =\int_{-\infty}^{\infty}
   \E^{-2t}\exp\{-\E^{-t}\}\,\D t
 =\int_0^\infty w\E^{-w}\,\D w=1.
\]
Thus $f$ is a probability density.  By Scheff\'e's theorem, see, e.g.,
\cite[Theorem~16.12]{Bi95}, the pointwise convergence
of the densities implies convergence in total variation, and hence
weak convergence.
\end{proof}

 The limiting density can be integrated explicitly, yielding a closed form for the joint distribution function.

\medskip

%%%%%%%%%%%%%%%%%%%%%%%%%%%%%%%%%%%%%%%%%%%%%%%%%%%%%%%%%%%%%%%%%%%%%%%%%%%%%%%%%%%%%%%%%%%%%%%%%%
%
% Proposition 2.2
%
\begin{proposition}[Joint distribution function]\label{prop:Daniels-joint-cdf}
%%%%%%%%%%%%%%%%%%%%%%%%%%%%%%%%%%%%%%%%%%%%%%%%%%%%%%%%%%%%%%%%%%%%%%%%%%%%%%%%%%%%%%%%%%%%%%%%%
For $t\in\mathbb R$ and $s\geq0$, put $w:=\E^{-t}$ and $A:=I_0(s)$.

Then
\begin{equation}\label{eq:joint-limit-cdf}
\begin{split}
 \PP(T\leq t,S\leq s)  ={}&(1+w)\E^{-w}  -\frac{1+Aw}{A^2}\E^{-Aw}.
\end{split}
\end{equation}
\end{proposition}
%%%%%%%%%%%%%%%%%%%%%%%%%%%%%%%%%%%%%%%%%%%%%%%%%%%%%%%%%%%%%%%%%%%%%%%%%%%%%%%%%%%%%%%%%%%%%%%%%%

\smallskip

\begin{proof}
Using $I_0'=I_1$, for fixed $u\in\mathbb R$ we have
\[
 \int_0^s I_1(v)  \exp\{-\E^{-u}I_0(v)\}\,\D v
 =\frac{\exp(-\E^{-u})- \exp(-\E^{-u}I_0(s))}{\E^{-u}}.
\]
Integrating this identity over $u\leq t$ and using $y=\E^{-u}$ gives
\[
\mathbb P(T\leq t,S\leq s) = \int_w^\infty y\{\E^{-y}-\E^{-Ay}\}\,\D y,
\]
where $w=\E^{-t}$ and $A=I_0(s)$. This is exactly \eqref{eq:joint-limit-cdf}.
\end{proof}

The marginal limits follow immediately from this representation.

\medskip

%%%%%%%%%%%%%%%%%%%%%%%%%%%%%%%%%%%%%%%%%%%%%%%%%%%%%%%%%%%%%%%%%%%%%%%%%%%%%%%%%%%%%%%%%%%%%%%
%
% Corollary 2.3
\begin{corollary}[Marginal limits]\label{cor:Daniels-marginals}
%
%%%%%%%%%%%%%%%%%%%%%%%%%%%%%%%%%%%%%%%%%%%%%%%%%%%%%%%%%%%%%%%%%%%%%%%%%%%%%%%%%%%%%%%%%%%%%%%%
The limiting radius variable $T$ has distribution function
\begin{equation}\label{eq:T-marginal}
 \PP(T\leq t)  =\exp\{-\E^{-t}\}(1+\E^{-t}),  \qquad t\in\mathbb R,
\end{equation}
and density
\[
 f_T(t)=\E^{-2t}\exp\{-\E^{-t}\}.
\]
Moreover,
\begin{equation}\label{eq:S-marginal}
 \PP(S\leq s)=1-\frac1{I_0(s)^2},  \qquad s\geq0, \end{equation}
and
\[
 f_S(s)=\frac{2I_1(s)}{I_0(s)^3},  \qquad s>0.
\]
\end{corollary}
%%%%%%%%%%%%%%%%%%%%%%%%%%%%%%%%%%%%%%%%%%%%%%%%%%%%%%%%%%%%%%%%%%%%%%%%%%%%%%%%%%%%%%%%%%%%%%%%

\smallskip

\begin{proof}
Letting $s\to\infty$ in \eqref{eq:joint-limit-cdf} gives \eqref{eq:T-marginal}.  Letting $t\to\infty$ gives
\eqref{eq:S-marginal}.  Differentiation yields the densities.
\end{proof}

\medskip

%%%%%%%%%%%%%%%%%%%%%%%%%%%%%%%%%%%%%%%%%%%%%%%%%%%%%%%%%%%%%%%%%%%%%%%%%%%%%%%%%%%%%%%%%%%%%%%%%
%
% Remark 2.4
\begin{remark}\label{rem:Daniels-centre}
%
%%%%%%%%%%%%%%%%%%%%%%%%%%%%%%%%%%%%%%%%%%%%%%%%%%%%%%%%%%%%%%%%%%%%%%%%%%%%%%%%%%%%%%%%%%%%%%%%%%%%%
The radius marginal \eqref{eq:T-marginal} is the non-Gumbel limit law that follows from Daniels' exact radius distribution.  Daniels also derived \eqref{eq:S-marginal} directly: in his notation the normalized
centre distance is $R=p\sqrt{2\log n}$, and his formula (5.3) gives the limiting distribution function $1-I_0(R)^{-2}$.  Theorem \ref{thm:Daniels-joint-limit} shows that these two marginal limits are
parts of a single explicit joint law.
\end{remark}

The joint limit admits a simple stochastic representation which is
useful both conceptually and for simulation.

\medskip

%%%%%%%%%%%%%%%%%%%%%%%%%%%%%%%%%%%%%%%%%%%%%%%%%%%%%%%%%%%%%%%%%%%%%%%%%%%%%%%%%%%%%%%%%%%%%%%%%%%
%
% Proposition 2.5
\begin{proposition}[Gamma--exponential representation]\label{prop:gamma-exp}
%
%%%%%%%%%%%%%%%%%%%%%%%%%%%%%%%%%%%%%%%%%%%%%%%%%%%%%%%%%%%%%%%%%%%%%%%%%%%%%%%%%%%%%%%%%%%%%%%%%%%
Let $W$ and $E$ be independent random variables with $W\sim\operatorname{Gamma}(2,1)$ and 
 $E\sim\operatorname{Exp}(1)$,
where the gamma distribution is parametrized by shape and rate.  Then
\[
%\begin{equation}\label{eq:gamma-exp-representation}
 T\edist-\log W,  \qquad  S\edist I_0^{-1}\left(1+\frac EW\right),
%\end{equation}
\]
where $I_0^{-1}$ denotes the inverse of $I_0:[0,\infty)\to[1,\infty)$. These two identities hold jointly, that is,
\[
 (T,S)\edist
 \left(-\log W,
 I_0^{-1}\left(1+\frac EW\right)\right).
\]
\end{proposition}
%%%%%%%%%%%%%%%%%%%%%%%%%%%%%%%%%%%%%%%%%%%%%%%%%%%%%%%%%%%%%%%%%%%%%%%%%%%%%%%%%%%%%%%%%%%%%%

\smallskip
\begin{proof}
Put $W:=\E^{-T}$  and $Y:=I_0(S)$.
The change of variables in \eqref{eq:joint-limit-density} gives
\[
 f_{W,Y}(w,y)=w^2\E^{-wy},
 \qquad w>0,\quad y>1.
\]
Now set $E=W(Y-1)$.  Since $y=1+E/W$, another change of variables gives
\[
 f_{W,E}(w,\xi)=w\E^{-w}\E^{-\xi},  \qquad w,\xi >0.
\]
This factors into the density of a $\operatorname{Gamma}(2,1)$ variable and that of an independent $\operatorname{Exp}(1)$ variable, which proves the assertion.
\end{proof}

\medskip

%%%%%%%%%%%%%%%%%%%%%%%%%%%%%%%%%%%%%%%%%%%%%%%%%%%%%%%%%%%%%%%%%%%%%%%%%%%%%%%%%%%%%%%%%%%%%%%%%
%
% Corollary 2.6
\begin{corollary}[Moments and covariance]
%
%%%%%%%%%%%%%%%%%%%%%%%%%%%%%%%%%%%%%%%%%%%%%%%%%%%%%%%%%%%%%%%%%%%%%%%%%%%%%%%%%%%%%%%%%%%%%%%%
Let $\gamma_{\rm E}$ denote Euler's constant, and put
\[
m_1:=\int_0^\infty \frac{ds}{I_0(s)^2},\qquad
m_2:=2\int_0^\infty \frac{s\,ds}{I_0(s)^2}, \qquad 
c_{TS}:=\int_0^\infty
\frac{\log I_0(s)}{I_0(s)^2}\,ds.
\]
Then $\BE T=\gamma_{\rm E}-1$, $\BE S=m_1$, and
\[
\operatorname{Cov} 
\begin{pmatrix}T\\ S\end{pmatrix}
= \begin{pmatrix}
\displaystyle \frac{\pi^2}{6}-1 & c_{TS}\\[1ex]
c_{TS} & m_2-m_1^2
\end{pmatrix}.
\]
Numerically, $m_1\approx1.36768$, $m_2-m_1^2\approx0.65499$, $c_{TS}\approx0.39696$,
and hence $\operatorname{Corr}(T,S)\approx0.61076$.
In particular, $T$ and $S$ are positively correlated and are not independent.
\end{corollary}
%%%%%%%%%%%%%%%%%%%%%%%%%%%%%%%%%%%%%%%%%%%%%%%%%%%%%%%%%%%%%%%%%%%%%%%%%%%%%%%%%%%%%%%%%%%%%%%%%%%%%%%%

\smallskip

\begin{proof} 
By Proposition~\ref{prop:gamma-exp}, $T=-\log W$ with
$W\sim\operatorname{Gamma}(2,1)$. Hence
$\mathbb E T=\gamma_{\rm E}-1$ and $\operatorname{Var}(T)=\pi^2/6-1$.
Moreover, Corollary~\ref{cor:Daniels-marginals}  gives
$\mathbb P(S>s)=1/I_0(s)^2$,
so that, by the tail integral formula,
$\mathbb E S=m_1$ and $\mathbb E S^2=m_2$. 
Finally, Proposition~\ref{prop:gamma-exp} yields, conditionally on $W=w$,
\[
\mathbb P(S>s\mid W=w)
=\exp\{-w(I_0(s)-1)\}.
\]
Using Fubini's theorem and the Gamma$(2,1)$ density,
\[
\mathbb E(TS)
= \int_0^\infty \frac{\log I_0(s)-\psi(2)}{I_0(s)^2}\, \D s.
\]
Since $\mathbb E T=-\psi(2)$, it follows that
\[
\operatorname{Cov}(T,S)
= \int_0^\infty \frac{\log I_0(s)}{I_0(s)^2}\, \D s =c_{TS}.
\]
The stated covariance matrix and numerical values follow.
\end{proof}

The preceding results concern the distance of the centre from the origin. Rotational invariance now yields the corresponding vector-valued limit.
Let $\Theta$ be uniformly distributed on $\mathbb S^1$ and independent
of $(T,S)$.

\medskip

%%%%%%%%%%%%%%%%%%%%%%%%%%%%%%%%%%%%%%%%%%%%%%%%%%%%%%%%%%%%%%%%%%%%%%%%%%%%%%%%%%%%%%%%%%%%%%%%%%%%%
%
% Corollary 2.6
%
\begin{corollary}[Vector-valued centre]\label{cor:vector-centre}
%%%%%%%%%%%%%%%%%%%%%%%%%%%%%%%%%%%%%%%%%%%%%%%%%%%%%%%%%%%%%%%%%%%%%%%%%%%%%%%%%%%%%%%%%%%%%%%%%%%%
As $n\to\infty$,
\begin{equation}\label{eq:vector-joint-limit}
 \left( a_n(R_n-a_n),\ a_nC_n \right) \vertk  \bigl(T,S\Theta\bigr).
\end{equation}
Equivalently, with $W,E,\Theta$ mutually independent,
\begin{equation}\label{eq:vector-representation}
 \bigl(T,S\Theta\bigr)  \edist  \left(  -\log W,  I_0^{-1}\left(1+\frac EW\right)\Theta  \right),
\end{equation}
where $W\sim\operatorname{Gamma}(2,1)$ and $E\sim\operatorname{Exp}(1)$.
\end{corollary}
%%%%%%%%%%%%%%%%%%%%%%%%%%%%%%%%%%%%%%%%%%%%%%%%%%%%%%%%%%%%%%%%%%%%%%%%%%%%%%%%%%%%%%%%%%%%%%%%%%%%%

\smallskip

\begin{proof}
For every $n$, rotational invariance of the Gaussian sample implies that, conditionally on $(R_n,\|C_n\|)$, the direction of $C_n$ is uniform on $\mathbb S^1$ and independent of these two rotation-invariant
quantities.  Combining this fact with Theorem~\ref{thm:Daniels-joint-limit} yields 
\eqref{eq:vector-joint-limit}.
The representation \eqref{eq:vector-representation} follows from
Proposition~\ref{prop:gamma-exp}.
\end{proof}

\medskip

%%%%%%%%%%%%%%%%%%%%%%%%%%%%%%%%%%%%%%%%%%%%%%%%%%%%%%%%%%%%%%%%%%%%%%%%%%%%%%%%%%%%%%%%%%%%%%
%
% Remark 2.7
\begin{remark}[Supporting configurations]\label{rem:supporting-configurations}
%
%%%%%%%%%%%%%%%%%%%%%%%%%%%%%%%%%%%%%%%%%%%%%%%%%%%%%%%%%%%%%%%%%%%%%%%%%%%%%%%%%%%%%%%%%%%%%%%%%%%
The two summands in Daniels' finite-sample density \eqref{eq:Daniels-joint-density} correspond to three-point and two-point support configurations, respectively. In the proof of Theorem~\ref{thm:Daniels-joint-limit}, the contribution of the second summand vanishes pointwise on the joint extreme-value scale, whereas
the contribution of the first converges pointwise to the density $f$ in \eqref{eq:joint-limit-density}. Since this contribution is nonnegative and $f$ has total mass one, Fatou's lemma implies that its
total mass tends to one. Hence the total mass of the second contribution tends to zero. Thus Daniels' formula also shows that the probability that the smallest enclosing circle is supported by two
points tends to zero, whereas the three-point configuration has probability tending to one.
\end{remark}

\medskip

\begin{remark}[Poisson interpretation]\label{rem:poisson-interpretation}
Theorem~\ref{thm:Daniels-joint-limit} will later identify explicitly
the value--argmin pair associated with the planar Gumbel Poisson
functional $G(z)$ defined in \eqref{eq:Hn-G}.  More precisely, once
one proves
\[
 \left(\frac{R_n-a_n}{b_n},\frac{C_n}{b_n}\right)
 \vertk (\Phi_G,Z_G)
\]
under the general Gumbel condition \eqref{eq:gumbel-tail}, comparison
with the Gaussian case gives
\[
 (\Phi_G,\|Z_G\|)\edist(T,S),
\]
and rotational invariance then yields the full vector law
$(\Phi_G,Z_G)\edist(T,S\Theta)$.  Hence the explicit Daniels law is not
restricted to the Gaussian model: it identifies the universal planar
Gumbel Poisson limit.
\end{remark}

%%%%%%%%%%%%%%%%%%%%%%%%%%%%%%%%%%%%%%%%%%%%%%%%%%%%%%%%%%%%%%%%%%%%%%%%%%%%%%%
%
% Section 3
\section{Joint asymptotics in the Gumbel regime}\label{sec:gumbel}
%
%%%%%%%%%%%%%%%%%%%%%%%%%%%%%%%%%%%%%%%%%%%%%%%%%%%%%%%%%%%%%%%%%%%%%%%%%%%%%%%

Assume throughout this section that there exist sequences $a_n>0$ and
$b_n>0$ such that
\begin{equation}\label{eq:gumbel-tail}
 n\PP(\rho>a_n+b_ns)\rightarrow \E^{-s},  \qquad s\in\mathbb R.
\end{equation}
As in Henze and Klatt~\cite{hekl}, condition~\eqref{eq:gumbel-tail}
implies $b_n/a_n\to0$. Recalling the representation $X_j=\rho_jU_j$ from Section~\ref{sec:intro}, put
\[
S_{n,j}:=\frac{\rho_j-a_n}{b_n},\qquad 1\le j\le n,
\]
and consider the marked point process
\[
 N_n:=\sum_{j=1}^n\varepsilon_{(S_{n,j},U_j)}.
\]
Then $N_n$ converges on bounded height strips to a Poisson point process
$N$ on $\mathbb R\times\mathbb S^{d-1}$ with intensity
$\E^{-s}\,\D s\,\sigma_{d-1}(\D u)$. Indeed, for $s<t$ and any Borel set $A\subset \mathbb{S}^{d-1}$,
\[
n\PP\{S_{n,1}\in(s,t],\,U_1\in A\}
\rightarrow
\bigl(\E^{-s}-\E^{-t}\bigr)\sigma_{d-1}(A)
= \int_s^t \E^{-x}\, \D x\,\sigma_{d-1}(A).
\]
Since the marked observations are independent, the standard
binomial-to-Poisson limit yields convergence of $N_n$, on every
bounded height strip, to a Poisson point process $N$ with the above intensity.
Define
\begin{equation}\label{eq:Hn-G}
 H_n(z):=\frac{1}{b_n}
 \left\{\max_{1\leq j\leq n}\|\rho_jU_j-b_nz\|-a_n\right\},
 \qquad  G(z):=\sup_{(s,u)\in N}(s-z\cdot u).
\end{equation}
The exact variational identities are
\begin{equation}\label{eq:Hn-value-argmin}
 \frac{R_n-a_n}{b_n}=\min_{z\in\mathbb R^d}H_n(z),
 \qquad  \frac{C_n}{b_n}=\argmin_{z\in\mathbb R^d}H_n(z).
\end{equation}

The joint limit will follow once the variational problem defined by $G$ is understood. We therefore begin with its basic geometric properties. In this connection, we call an atom $(s,u)\in N$ active at $z$ if
$s-z\cdot u=G(z)$.

\medskip

%%%%%%%%%%%%%%%%%%%%%%%%%%%%%%%%%%%%%%%%%%%%%%%%%%%%%%%%%%%%%%%%%%%%%%%%%%%%%%%%%%%%%%%%%%
%
% Lemma 3.1
\begin{lemma}[Geometry of the Gumbel limit]\label{lem:gumbel-geometry}
%
%%%%%%%%%%%%%%%%%%%%%%%%%%%%%%%%%%%%%%%%%%%%%%%%%%%%%%%%%%%%%%%%%%%%%%%%%%%%%%%%%%%%%%%%%%%%%
Almost surely, $G$ is real-valued, $1$-Lipschitz, convex,  and $G(z) \to \infty$ as $\|z\| \to \infty$.  It has a unique minimizer $Z_G$; writing
\[
 \Phi_G:=G(Z_G),
\]
exactly $d+1$ atoms of $N$ are active at $Z_G$.  Their directions are affinely independent and contain the origin in the interior of their convex hull.  In particular, there is an almost surely positive random
variable $\eta_G$ such that 
\begin{equation}\label{eq:gumbel-linear-growth}
 G(z)\geq \Phi_G+\eta_G\|z-Z_G\|,  \qquad z\in\mathbb R^d.
\end{equation}
\end{lemma}
%%%%%%%%%%%%%%%%%%%%%%%%%%%%%%%%%%%%%%%%%%%%%%%%%%%%%%%%%%%%%%%%%%%%%%%%%%%%%%%%%%%%%%%%%%%%

\smallskip

\begin{proof}
Let $S_*:=\sup_{(s,u)\in N} s$. Since
$\PP(S_*\le t) = \PP\{N((t,\infty)\times \mathbb{S}^{d-1})=0\} = \exp\{-\E^{-t}\}$,
$t\in\mathbb R$, $S_*$ has the standard Gumbel distribution and is therefore finite
almost surely. On the ball $\{\|z\|\le M\}$, atoms $(s,u)$ of $N$ with
$s<S_*-2M$ cannot affect the supremum in the definition of $G$.
 Only finitely many atoms remain above this level.  Hence $G$ is
real-valued and, locally, the maximum of finitely many affine functions. The $1$-Lipschitz property and convexity follow immediately.

Let $v_0,\ldots,v_d\in \mathbb{S}^{d-1}$ be the vertices of a regular
simplex centred at the origin. Then, for some $\delta_0>0$,  
$\delta_0 \mathbb{B}^d\subset\operatorname{conv}\{v_0,\ldots,v_d\}$.
By continuity, there exist pairwise disjoint open spherical
neighbourhoods $C_i$ of $v_i$, $0\le i\le d$, and a constant
$\delta>0$ such that
$\delta \mathbb{B}^d\subset\operatorname{conv}\{u_0,\ldots,u_d\}$
for every choice $u_i\in C_i$.
Since $\sigma_{d-1}(C_i)>0$ and
$\int_{\mathbb R}\E^{-s}\, \D s=\infty$,
each $C_i$ contains almost surely the direction of at least one atom
$(s_i,u_i)$ of $N$. For every $z\in\mathbb R^d$,
$\max_{0\le i\le d}(-z\cdot u_i)\ge \delta\|z\|$,
and therefore
\[
G(z) \ge \max_{0\le i\le d}(s_i-z\cdot u_i) \ge \min_{0\le i\le d}s_i+\delta\|z\|.
\]
Thus $G$ attains its minimum, and $G(z) \to \infty$ as $\|z\| \to \infty$.

We next record the relevant general-position properties.  Almost surely, every collection of at most $d$ distinct atom directions is linearly independent, every collection of $d+1$ distinct atom
directions is affinely independent, and no $d+2$ affine functions $s-z\cdot u$ associated with distinct atoms take the same value at a common point $z$.  The first two assertions follow from absolute
continuity of spherical measure.  For the third, the first $d+1$ equalities determine $(z,v)\in\mathbb R^{d+1}$ uniquely whenever the directions are affinely independent; conditionally on these atoms, the
height of any further atom has a continuous distribution.  Applying this argument on bounded height strips and then taking a countable union proves the claim.

Fix a realization of $N$ for which the above general-position properties
hold and for which $G$ attains its minimum, and let $z$ be a minimizer
of $G$. In a neighbourhood of $z$, only finitely many affine functions
can be relevant.
If zero did not belong to the convex hull of the directions of the
active functions, the strict separation theorem from convex geometry
would yield a vector $h\in\mathbb R^d$ and a constant $c>0$ such that
$h\cdot u\ge c$
for every active direction $u$. Hence, for every active atom $(s,u)$
and every $\varepsilon>0$,
\[
s-(z+\varepsilon h)\cdot u
=
s-z\cdot u-\varepsilon h\cdot u
\le
G(z)-\varepsilon c.
\]
Among the finitely many affine functions relevant in a neighbourhood
of $z$, those that are not active at $z$ have a strictly positive gap
below $G(z)$. Hence, for sufficiently small $\varepsilon>0$, all
relevant affine functions have value strictly smaller than $G(z)$ at
$z+\varepsilon h$, contradicting the minimality of $z$.
Thus
$0\in\operatorname{conv}\{u:(s,u)\text{ is active at }z\}$.

By linear independence, at least $d+1$ atoms are active, whereas the general-position property excludes $d+2$ active atoms.  Thus exactly $d+1$ atoms are active.  Their directions are affinely independent;
moreover, the origin lies in the interior of their convex hull, since otherwise it would belong to the convex hull of at most $d$ of them.

Let $(s_j,u_j)$, $0\leq j\leq d$, be the active atoms.  Since $0\in\operatorname{int}\operatorname{conv}\{u_0,\ldots,u_d\}$, there
is $\eta>0$ such that $\eta \mathbb{B}^d\subset  \operatorname{conv}\{-u_0,\ldots,-u_d\}$.
For every $h\in\mathbb R^d$,
\[
 G(z+h)\geq \max_{0\leq j\leq d}  \{s_j-(z+h)\cdot u_j\}  =G(z)+\max_{0\leq j\leq d}(-h\cdot u_j)
 \geq G(z)+\eta\|h\|.
\]
This proves both uniqueness and~\eqref{eq:gumbel-linear-growth}.
\end{proof}

To pass from local uniform convergence of the objectives to convergence of their minimizers, we use the following elementary localization principle.

\medskip

%%%%%%%%%%%%%%%%%%%%%%%%%%%%%%%%%%%%%%%%%%%%%%%%%%%%%%%%%%%%%%%%%%%%%%%%%%%%%%%%%%%%%%%%%%%%%%%%%%%
%
% Lemma 3.2
\begin{lemma}[Localization by convexity]\label{lem:convex-localization}
%
%%%%%%%%%%%%%%%%%%%%%%%%%%%%%%%%%%%%%%%%%%%%%%%%%%%%%%%%%%%%%%%%%%%%%%%%%%%%%%%%%%%%%%%%%%%%%%%%
Let $f_n$ be real-valued convex functions on $\mathbb R^d$, each attaining
its minimum, and suppose that $f_n\to f$ uniformly on compact sets.
Assume that $f(z)\rightarrow\infty$ as $\|z\|\rightarrow\infty$,
and that $f$ has a unique minimizer $z_*$. Then every choice
$z_n\in\operatorname{argmin} f_n$ satisfies
\[
\bigl(\min f_n,z_n\bigr)\rightarrow
\bigl(\min f,z_*\bigr).
\]
If, in addition, for some $\eta >0$,
\[
 f(z)\geq f(z_*)+\eta\|z-z_*\|,
\]
then, on every ball containing $z_*$ and all sufficiently large $z_n$,
\[
 \|z_n-z_*\|  \leq \frac{2}{\eta}\sup|f_n-f|,  \qquad  |\min f_n-\min f|\leq\sup|f_n-f|,
\]
where the suprema are over that ball.
\end{lemma}
%%%%%%%%%%%%%%%%%%%%%%%%%%%%%%%%%%%%%%%%%%%%%%%%%%%%%%%%%%%%%%%%%%%%%%%%%%%%%%%%%%%%%%%%%%%%%%%

\smallskip

\begin{proof}
Choose $M$ so large that $f(z)>f(0)+1$ whenever $\|z\|=M$. For all sufficiently large $n$, local uniform convergence gives $f_n(z)>f_n(0)$ on this sphere.  A minimizer of $f_n$ outside 
the closed ball $B:= \{z\in\mathbb R^d:\|z\|\le M\}$
would, by convexity along the segment joining it to the origin, produce a point on the sphere with value at most $f_n(0)$, a contradiction. Thus, for all sufficiently large $n$, every minimizer of $f_n$ lies in $B$. Uniform convergence there and uniqueness of $z_*$ identify every subsequential
limit and yield convergence of the minimum values.

For the quantitative estimate, put $\varepsilon_n:=\sup_{z\in B}|f_n(z)-f(z)|$.
From
$f_n(z_n)\leq f_n(z_*)$ we obtain $f(z_n)\leq f(z_*)+2\varepsilon_n$ and hence
$\eta\|z_n-z_*\|\leq2\varepsilon_n$.  The bound for the minimum values is immediate.
\end{proof}

\smallskip

Recalling $H_n$ and $G$ from \eqref{eq:Hn-G}, we can now combine the point-process convergence with this localization principle to obtain the joint Gumbel limit.

\medskip

%%%%%%%%%%%%%%%%%%%%%%%%%%%%%%%%%%%%%%%%%%%%%%%%%%%%%%%%%%%%%%%%%%%%%%%%%%%%%%%%%%%%%%%%%%%%%%%%%%%
%
% Theorem 3.3
\begin{theorem}[Joint Gumbel limit]\label{thm:gumbel-joint}
%
%%%%%%%%%%%%%%%%%%%%%%%%%%%%%%%%%%%%%%%%%%%%%%%%%%%%%%%%%%%%%%%%%%%%%%%%%%%%%%%%%%%%%%%%%%%%%%%%%%%%
Under~\eqref{eq:gumbel-tail},
\[%begin{equation}\label{eq:gumbel-function-convergence}
 H_n\vertk G  \qquad\text{in }C_{\rm loc}(\mathbb R^d),
\]%end{equation}
where $C_{\rm loc}(\mathbb R^d)$ denotes the space of continuous
real-valued functions on $\mathbb R^d$, endowed with the topology of
uniform convergence on compact sets. Consequently, 
\begin{equation}\label{eq:gumbel-joint-limit}
 \left( \frac{R_n-a_n}{b_n},\frac{C_n}{b_n} \right)  \vertk(\Phi_G,Z_G).
\end{equation}
\end{theorem}
%%%%%%%%%%%%%%%%%%%%%%%%%%%%%%%%%%%%%%%%%%%%%%%%%%%%%%%%%%%%%%%%%%%%%%%%%%%%%%%%%%%%%%%%%%%%%%%%%

\smallskip

\begin{proof}
Fix $M>0$ and write $B_M:=\{z:\|z\|\leq M\}$. For $L<A$, let $N_n^{L,A}$ and $N^{L,A}$ denote the
restrictions of $N_n$ and $N$ to $[L,A]\times\mathbb S^{d-1}$.  Put
\[
 \psi_n(s,u;z)
 :=
 \frac{\|(a_n+b_ns)u-b_nz\|-a_n}{b_n}
\]
and define, for $z\in B_M$,
\[
 \widetilde H_{n,L,A}(z)  :=  \max\left\{  L+M,\,
 \max_{\substack{1\leq j\leq n\\ L\leq S_{n,j}\leq A}}  \psi_n(S_{n,j},U_j;z)  \right\},
\]
as well as
\[
 \widetilde G_{L,A}(z)  :=  \max\left\{  L+M,\,
 \max_{\substack{(s,u)\in N\\ L\leq s\leq A}}  (s-z\cdot u)  \right\}.
\]
The inner maximum is interpreted as $-\infty$ if the corresponding
point configuration is empty.  Thus both functions are well defined
as elements of $C(B_M)$, where $C(B_M)$ denotes the space of continuous real-valued functions
on $B_M$, equipped with the supremum norm. The point-process convergence implied by 
\eqref{eq:gumbel-tail} gives $ N_n^{L,A}\vertk N^{L,A}$ whenever the boundary of the strip contains no atom of $N$, which is an event of probability one.  
The limiting restricted process is finite and simple.   By
Proposition~3.13 of Resnick~\cite{resn}, at every such limiting
configuration the atoms of a vaguely convergent sequence of point
measures can, for all sufficiently large indices, be labelled so that
the corresponding heights and directions converge. Together with the
continuous mapping theorem, this implies the required convergence of
the corresponding finite configurations.

Moreover, uniformly for $s\in[L,A]$, $u\in\mathbb S^{d-1}$ and
$z\in B_M$, we have
$\psi_n(s,u;z)=s-z\cdot u+o(1)$,
where the $o(1)$ is deterministic. Consequently,
$\widetilde H_{n,L,A}\vertk \widetilde G_{L,A}$ in $C(B_M)$.

Let $S_{n,*}:=\max_{1\le j\le n}S_{n,j}$. On the event
$\{L+2M<S_{n,*}\le A\}$, we have
$H_n(z)=\widetilde H_{n,L,A}(z)$ for $z\in B_M$.
Indeed, every index $j$ with $S_{n,j}<L$ contributes at most $L+M$,
whereas an index $j_*$ satisfying $S_{n,j_*}=S_{n,*}$ contributes at
least
$S_{n,*}-M>L+M$.
The condition $S_{n,*}\le A$ excludes terms with height above the
truncation level. The same argument applies to the limiting process.
Furthermore, $\PP(S_{n,*}\leq L+2M)  \rightarrow  \exp\{-\E^{-(L+2M)}\}$, whereas
\[
 \limsup_{n\to\infty}\PP(S_{n,*}>A)\leq \E^{-A}.
\]
The corresponding bounds hold for the limiting Poisson process. Letting first $L\to-\infty$ and then $A\to\infty$, the converging-together theorem
(see, e.g., \cite[Theorem~3.2]{bi})
therefore yields
$H_n|_{B_M}\vertk G|_{B_M}$ in $C(B_M)$.
Since the convergence holds in $C(B_M)$ for every positive integer $M$, it follows that
$H_n \vertk G$ in $C_{\mathrm{loc}}(\mathbb R^d)$.
Each $H_n$ is real-valued and convex, satisfies
$H_n(z)\to\infty$ as $\|z\|\to\infty$,
and has the unique minimizer $C_n/b_n$. Lemma~\ref{lem:gumbel-geometry} shows that, almost surely,
$G(z)\to\infty$ as $\|z\|\to\infty$,
and that $G$ has a unique minimizer. Thus Lemma~\ref{lem:convex-localization}, or equivalently the continuous
mapping theorem applied to the value--argmin map at $G$, together with
\eqref{eq:Hn-value-argmin}, proves~\eqref{eq:gumbel-joint-limit}.
\end{proof}

We next record a consequence for the supporting observations.  Let
$K_n$ denote the number of sample points lying on the boundary of the
smallest enclosing ball.

\medskip

%%%%%%%%%%%%%%%%%%%%%%%%%%%%%%%%%%%%%%%%%%%%%%%%%%%%%%%%%%%%%%%%%%%%%%%%%%%%%%%%%%%%%%%%%%%%%%%
%
% Lemma 3.4
\begin{lemma}[Stability of a nondegenerate active simplex]\label{lem:active-stability}
%
%%%%%%%%%%%%%%%%%%%%%%%%%%%%%%%%%%%%%%%%%%%%%%%%%%%%%%%%%%%%%%%%%%%%%%%%%%%%%%%%%%%%%%%%%%%%
Let $m \ge d+1$ be an integer, and let
\[
\ell_j(z)=\alpha_j+u_j\cdot z,
\qquad 1\le j\le m,
\]
where $\alpha_j\in\mathbb R$ and $u_j\in\mathbb R^d$. Put
$F(z)=\max_j \ell_j(z)$.
Suppose that $F$ has a minimizer $z^\ast$ and that its active set $I$
has cardinality $d+1$. Assume that the vectors $(u_j)_{j\in I}$
are affinely independent and
$0\in \operatorname{int}\operatorname{conv}\{u_j:j\in I\}$.

Let $f_{n,j}$ be convex functions which are $C^1$ in a neighbourhood of $z_*$ and satisfy $f_{n,j}\to\ell_j$ in $C^1$ there for every $j$. Then, for all sufficiently large $n$, there is a unique point $z_n$ near $z_*$ at which precisely the constraints indexed by $I$ are
active, and $z_n$ is the unique global minimizer of $\max_j f_{n,j}$. Moreover, $z_n\to z_*$ and the corresponding minimum values converge.
\end{lemma}
%%%%%%%%%%%%%%%%%%%%%%%%%%%%%%%%%%%%%%%%%%%%%%%%%%%%%%%%%%%%%%%%%%%%%%%%%%%%%%%%%%%%%%%%%%%%%%

\smallskip

\begin{proof}
Relabel $I=\{1,\ldots,d+1\}$ and introduce the additional variable $v$.
The limiting system
\[
\ell_j(z)-v=0,\qquad 1\le j\le d+1,
\]
has the solution $(z^\ast,F(z^\ast))$. Its Jacobian at
$(z^\ast,F(z^\ast))$ has rows
$(u_j^\top,-1)$ for $1\le j\le d+1$.
This matrix is nonsingular precisely because the vectors
$u_1,\ldots,u_{d+1}$ are affinely independent.
 The inverse function theorem and the
$C^1$ convergence therefore yield a unique nearby solution $(z_n,v_n)$ with $(z_n,v_n)\to(z_*,F(z_*))$.

Since the inactive set is finite, the inactive affine functions have a strictly positive gap at $z_*$.  Local uniform convergence preserves this gap at $z_n$ for all sufficiently large $n$.  Furthermore,
$0\in\operatorname{int}\operatorname{conv}\{u_j:j\in I\}$ is stable under small perturbations.  Hence, for some $\eta>0$ and all large $n$,
$\eta \mathbb{B}^d\subset  \operatorname{conv}\{\nabla f_{n,j}(z_n):j\in I\}$.
For $h\in\mathbb R^d$, convexity gives
\[
 \max_j f_{n,j}(z_n+h)  \geq \max_{j\in I}  \{v_n+\nabla f_{n,j}(z_n)\cdot h\}
 \geq v_n+\eta\|h\|.
\]
Thus $z_n$ is the unique global minimizer, and the strict inactive gaps
show that its active set is exactly $I$.
\end{proof}

\smallskip

We can now transfer the nondegenerate active-simplex structure of the
Poisson limit to the supporting observations of the sample.

\medskip

%%%%%%%%%%%%%%%%%%%%%%%%%%%%%%%%%%%%%%%%%%%%%%%%%%%%%%%%%%%%%%%%%%%%%%%%%%%%%%%%%%%%%%%%%%%%
%
% Proposition 3.5
\begin{proposition}[Asymptotic supporting simplex]\label{prop:gumbel-support}
%
%%%%%%%%%%%%%%%%%%%%%%%%%%%%%%%%%%%%%%%%%%%%%%%%%%%%%%%%%%%%%%%%%%%%%%%%%%%%%%%%%%%%%%%%%
Under~\eqref{eq:gumbel-tail}, we have
\[
\lim_{n\to \infty}  \PP(K_n=d+1)\ =  \ 1.
\]
More precisely, define the point measure of active observations by
\[
 \mathcal A_n:=  \sum_{\substack{1\leq j\leq n\\ \|X_j-C_n\|=R_n}}
 \varepsilon_{(S_{n,j},U_j)}
\]
and the limiting active point measure by
\[
 \mathcal A:=  \sum_{\substack{(s,u)\in N\\ s-Z_G\cdot u=\Phi_G}}  \varepsilon_{(s,u)}.
\]
Then, with the vague topology on finite point measures,
\[
 \left( \frac{R_n-a_n}{b_n},\frac{C_n}{b_n},\mathcal A_n \right)
 \vertk  (\Phi_G,Z_G,\mathcal A).
\]
Almost surely, $\mathcal A$ has exactly $d+1$ atoms.
\end{proposition}
%%%%%%%%%%%%%%%%%%%%%%%%%%%%%%%%%%%%%%%%%%%%%%%%%%%%%%%%%%%%%%%%%%%%%%%%%%%%%%%%%%%%%%%%%%%%%%%%

\smallskip

\begin{proof}
Choose $M>0$ and a compact height strip $[L,A]$, as in the proof of
Theorem~3.3, and put $B_M:=\{z\in\mathbb R^d:\|z\|\le M\}$.
Choose $M$, $L$, and $A$ so that, with probability arbitrarily close to one,
$Z_G\in B_M$, the restriction of $G$ to $B_M$ is determined by the atoms
with heights in $[L,A]$, and, asymptotically, the restriction of $H_n$ to
$B_M$ is determined by the corresponding sample points with heights in
$[L,A]$.

By the Skorokhod representation theorem
(see, e.g., Billingsley~\cite[Theorem 6.7]{bi}),
we may work on a probability space on which the restricted finite
point measures converge almost surely.
 By Lemma~\ref{lem:gumbel-geometry}, the limiting active set consists of exactly $d+1$ atoms, their directions are affinely independent, and the origin lies in the interior of their convex hull.  All other retained
constraints have a strictly positive gap at $Z_G$.

For an atom $(S_{n,j},U_{n,j})$ corresponding to the limiting atom
$(s_j,u_j)$, define
\[
 F_{n,j}(z):=
 \frac{\|(a_n+b_nS_{n,j})U_{n,j}-b_nz\|-a_n}{b_n}.
\]
On a neighbourhood of $Z_G$, these functions converge in $C^1$ to $s_j-z\cdot u_j$.  Indeed, after factoring out $a_n$, the function $F_{n,j}(z)$ can be
written as a smooth function of $b_n/a_n$, $S_{n,j}$, $U_{n,j}$,
and $z$.    Applying Lemma~\ref{lem:active-stability} to the retained constraints gives, for all sufficiently large $n$, a unique nearby point $z_n$ at which exactly the $d+1$ corresponding 
constraints are active.  Their maximum has a strict global minimum at $z_n$.

The strict gaps exclude all other retained constraints.  The truncation inequalities from the proof of
Theorem~\ref{thm:gumbel-joint} exclude the discarded constraints at $z_n$. Hence the function $H_n$ agrees at $z_n$ with the maximum of the selected
$d+1$ constraints, while everywhere it is at least this maximum.  Thus $z_n$ is the unique global minimizer of $H_n$ and exactly these $d+1$ sample observations support the smallest enclosing ball.  
This correspondence also gives convergence of their heights and original directions.  
Since, on the event that $Z_G\in B_M$ and that both $G$ and $H_n$
on $B_M$ are determined by the atoms and sample points, respectively,
with heights in $[L,A]$, the active atoms are exactly those corresponding
to the limiting active set, their point measure converges vaguely to
$\mathcal A$ jointly with the value and the optimizer.
Letting the probability of the localization complement tend to zero proves the assertion.
\end{proof}

\smallskip

Note that the first assertion of Proposition~\ref{prop:gumbel-support} answers Problem 7.6 of Henze and Klatt~\cite{hekl}. In fact, it gives the stronger joint convergence of the active normalized observations to the active Poisson simplex.

In dimension two the limiting law can now be identified explicitly by
Section~\ref{sec:daniels}.

\medskip

%%%%%%%%%%%%%%%%%%%%%%%%%%%%%%%%%%%%%%%%%%%%%%%%%%%%%%%%%%%%%%%%%%%%%%%%%%%%%%%%%%%%%%%%%%%%%%%%
%
% Corollary 3.6
\begin{corollary}[Planar universality]\label{cor:gumbel-planar}
%
%%%%%%%%%%%%%%%%%%%%%%%%%%%%%%%%%%%%%%%%%%%%%%%%%%%%%%%%%%%%%%%%%%%%%%%%%%%%%%%%%%%%%%%%%%%%%%
If $d=2$ and~\eqref{eq:gumbel-tail} holds, then
\[
 \left( \frac{R_n-a_n}{b_n},\frac{C_n}{b_n} \right)  \vertk (T,S\Theta),
\]
where $\Theta$ is uniform on $\mathbb S^1$ and independent of $(T,S)$, and $(T,S)$ has density
\[
 f(t,s)=\E^{-3t}I_1(s)  \exp\{-\E^{-t}I_0(s)\},  \qquad t\in\mathbb R,\ s>0.
\]
\end{corollary}
%%%%%%%%%%%%%%%%%%%%%%%%%%%%%%%%%%%%%%%%%%%%%%%%%%%%%%%%%%%%%%%%%%%%%%%%%%%%%%%%%%%%%%%%%%%%%%%%%%

\smallskip

\begin{proof}
For the standard planar Gaussian distribution one may take $a_n=\sqrt{2\log n}$ and $b_n=1/a_n$.  Theorem
\ref{thm:Daniels-joint-limit} and Corollary~\ref{cor:vector-centre}
therefore identify the law of the universal Poisson value--argmin pair
$(\Phi_G,Z_G)$ as $(T,S\Theta)$.  Theorem~\ref{thm:gumbel-joint} then
gives the assertion for every radial distribution satisfying
\eqref{eq:gumbel-tail}.
\end{proof}

%%%%%%%%%%%%%%%%%%%%%%%%%%%%%%%%%%%%%%%%%%%%%%%%%%%%%%%%%%%%%%%%%%%%%%%%%%%%%%%
%
% Section 4
\section{Joint asymptotics in the Fr\'echet regime}\label{sec:frechet}
%
%%%%%%%%%%%%%%%%%%%%%%%%%%%%%%%%%%%%%%%%%%%%%%%%%%%%%%%%%%%%%%%%%%%%%%%%%%%%%%%

We now assume that the radial distribution has a regularly varying
tail with index $-\alpha$, $\alpha>0$:
\begin{equation}\label{eq:frechet-tail}
 \PP(\rho>s)=s^{-\alpha}L(s),
 \qquad s\to\infty,
\end{equation}
where $L$ is slowly varying.  Put
\[
%\begin{equation}\label{eq:frechet-cn}
 c_n:=\inf\{s>0:\PP(\rho\leq s)\geq1-1/n\}.
%\end{equation}
\]
Then $c_n\to\infty$ and $n\PP(\rho>c_n)\to1$. As shown in Section~4 of Henze and Klatt~\cite{hekl},
\[
 \sum_{j=1}^n\varepsilon_{c_n^{-1}X_j}  \vertk  \sum_{k\geq1}\varepsilon_{P_k}
\]
vaguely on $\mathbb R^d\setminus\{0\}$, where $\{P_k:k\ge 1\}$ is the set of atoms of the limiting Poisson
point process with intensity measure $\nu$ given by
\[
%\begin{equation}\label{eq:frechet-intensity}
 \nu\left\{x:\|x\|>r,\ \frac{x}{\|x\|}\in A\right\}  =r^{-\alpha}\sigma_{d-1}(A),  \qquad r>0.
%\end{equation}
\]
Let $ \mathcal K:=\{0\}\cup\{P_k:k\geq1\}$.
The set $\mathcal K$ is almost surely compact: the Poisson cloud is bounded, has finitely many points outside every neighbourhood of the origin, and can accumulate only at the origin.
For a nonempty compact set $A\subset\mathbb R^d$, write
\[
 r(A):=\min_{y\in\mathbb R^d}\max_{x\in A}\|x-y\|
\]
for its circumradius and $c(A)$ for the uniquely determined centre of
its smallest enclosing ball.

To transfer point-process convergence to the radius--centre pair, we first establish continuity of the smallest-enclosing-ball map under Hausdorff perturbations.

\medskip
%%%%%%%%%%%%%%%%%%%%%%%%%%%%%%%%%%%%%%%%%%%%%%%%%%%%%%%%%%%%%%%%%%%%%%%%%%%%%%%%%%%%%%%%%%%%%%%
%
% Lemma 4.1
\begin{lemma}[Hausdorff stability of radius and centre]\label{lem:hausdorff-stability}
%
%%%%%%%%%%%%%%%%%%%%%%%%%%%%%%%%%%%%%%%%%%%%%%%%%%%%%%%%%%%%%%%%%%%%%%%%%%%%%%%%%%%%%%%%%%%%%%%%%%%%%%

Let $A,B$ be nonempty compact subsets of $\mathbb R^d$, and let $d_H(A,B)=:h $ denote their Hausdorff distance.  Then
\begin{equation}\label{eq:radius-stability}
 |r(A)-r(B)|\leq h
\end{equation}
and
\begin{equation}\label{eq:centre-stability}
 \|c(A)-c(B)\|^2  \leq4r(A)h+4h^2.
\end{equation}
Consequently, the map $A\longmapsto(r(A),c(A))$ is continuous on the space of nonempty compact subsets of
$\mathbb R^d$, equipped with the Hausdorff metric.
\end{lemma}
%%%%%%%%%%%%%%%%%%%%%%%%%%%%%%%%%%%%%%%%%%%%%%%%%%%%%%%%%%%%%%%%%%%%%%%%%%%%%%%%%%%%%%%%%%%%%%%%%%%%%%%%

\smallskip

\begin{proof}
For $A$ fixed, put $F_A(y):=\max_{x\in A}\|x-y\|^2$,
 $a:=c(A)$, and $r:=r(A)$. We first show that
\begin{equation}\label{eq:strong-centre-ineq}
 F_A(y)\geq r^2+\|y-a\|^2,  \qquad y\in\mathbb R^d.
\end{equation}
For $0<\lambda<1$, let
$a_\lambda=(1-\lambda)a+\lambda y$.  The identity
\[
 \|x-a_\lambda\|^2  =(1-\lambda)\|x-a\|^2  +\lambda\|x-y\|^2  -\lambda(1-\lambda)\|y-a\|^2
\]
shows, after taking maxima over $x\in A$, that
\[
 F_A(a_\lambda)  \leq(1-\lambda)F_A(a)+\lambda F_A(y)  -\lambda(1-\lambda)\|y-a\|^2.
\]
Since $a$ minimizes $F_A$, we have $F_A(a)\leq F_A(a_\lambda)$. Dividing by $\lambda$ and letting $\lambda\downarrow0$ proves
\eqref{eq:strong-centre-ineq}.

By the definition of Hausdorff distance, for every $y$,
\[
 \max_{x\in A}\|x-y\|  \leq \max_{x\in B}\|x-y\|+h.
\]
Taking minima over $y$ and interchanging $A$ and $B$ proves
\eqref{eq:radius-stability}.  Taking $y=c(B)$ in \eqref{eq:strong-centre-ineq} gives
\[
 r(A)^2+\|c(A)-c(B)\|^2  \leq F_A(c(B))  \leq(r(B)+h)^2  \leq(r(A)+2h)^2,
\]
which is~\eqref{eq:centre-stability}.  Continuity follows immediately.
\end{proof}

\smallskip
We next upgrade the vague convergence of the extremal point process to Hausdorff convergence of the associated compact sets.

\medskip

%%%%%%%%%%%%%%%%%%%%%%%%%%%%%%%%%%%%%%%%%%%%%%%%%%%%%%%%%%%%%%%%%%%%%%%%%%%%%%%%%%%%%%%%%%%%%%%%%%%%%
%
% Lemma 4.2
%
\begin{lemma}[Compact-set convergence]\label{lem:frechet-compact}
%
%%%%%%%%%%%%%%%%%%%%%%%%%%%%%%%%%%%%%%%%%%%%%%%%%%%%%%%%%%%%%%%%%%%%%%%%%%%%%%%%%%%%%%%%%%%%%%%%%%%
With $\mathcal K_n:=\{0,c_n^{-1}X_1,\ldots,c_n^{-1}X_n\}$,
we have
\[
 \mathcal K_n\vertk\mathcal K
\]
in the space of nonempty compact subsets of $\mathbb R^d$, equipped
with the Hausdorff metric.  The same convergence holds if the origin is
omitted from $\mathcal K_n$.
\end{lemma}
%%%%%%%%%%%%%%%%%%%%%%%%%%%%%%%%%%%%%%%%%%%%%%%%%%%%%%%%%%%%%%%%%%%%%%%%%%%%%%%%%%%%%%%%%%%%%%%

\smallskip

\begin{proof}
Fix $0<\varepsilon<L<\infty$. On the annulus $\{x:\varepsilon\le \|x\|\le L\}$, vague
point-process convergence and Proposition~3.13 of
Resnick~\cite{resn} imply that, whenever no limiting atom lies
on either boundary sphere, the atoms can be labelled so that the
corresponding atom locations converge.
 Hence the associated finite point configurations converge in
the Hausdorff metric. This holds almost surely for fixed $\varepsilon,L$. All points inside $\varepsilon \mathbb{B}^d$ are within distance $\varepsilon$ of the adjoined origin.  Moreover,
\[
 \PP\left(\max_{1\leq j\leq n}\frac{\rho_j}{c_n}>L\right)
 \leq n\PP(\rho>Lc_n)\rightarrow L^{-\alpha},
\]
whereas
\[
 \PP\left(\sup_{k\geq1}\|P_k\|>L\right)  =1-\exp\{-L^{-\alpha}\}\rightarrow0  \qquad(L\to\infty).
\]
For $0<\varepsilon<L$, define
\[
 \mathcal K_n^{\varepsilon,L}
 :=  \{0\}\cup  \{c_n^{-1}X_j:\varepsilon\leq \|c_n^{-1}X_j\|\leq L\}, \qquad
 \mathcal K^{\varepsilon,L}
 :=  \{0\}\cup  \{P_k:\varepsilon\leq \|P_k\|\leq L\}.
\]
The preceding finite-configuration convergence implies
$\mathcal K_n^{\varepsilon,L}  \vertk \mathcal K^{\varepsilon,L}$
in the Hausdorff metric. On the event $ \max_{1\leq j\leq n} \rho_j/c_n\leq L$,
all points omitted from $\mathcal K_n^{\varepsilon,L}$ lie within
distance $\varepsilon$ of the origin, and hence
$d_H(\mathcal K_n,\mathcal K_n^{\varepsilon,L})\leq\varepsilon$.
Likewise, on the event $\sup_k\|P_k\|\leq L$, $d_H(\mathcal K,\mathcal K^{\varepsilon,L})\leq\varepsilon$.
Together with
\[
 \limsup_{n\to\infty}  \PP\left(  \max_{1\leq j\leq n}\frac{\rho_j}{c_n}>L  \right)  \leq L^{-\alpha}
\]
and
\[
 \PP\left(\sup_{k\geq1}\|P_k\|>L\right)  =1-\exp\{-L^{-\alpha}\},
\]
a converging-together argument, first letting $n\to\infty$ and then $\varepsilon\downarrow0$ and $L\to\infty$, proves the asserted Hausdorff convergence.
Finally, if
$\mathcal S_n:=\{c_n^{-1}X_1,\ldots,c_n^{-1}X_n\}$, then
\[
 d_H(\mathcal S_n,\mathcal K_n)  \leq \frac{\|X_1\|}{c_n}\stk0,
\]
so adjoining the origin does not affect the limit.
\end{proof}

\smallskip
The joint Fréchet limit is now an immediate consequence of Lemmas~\ref{lem:hausdorff-stability}
and~\ref{lem:frechet-compact}.

\medskip

%%%%%%%%%%%%%%%%%%%%%%%%%%%%%%%%%%%%%%%%%%%%%%%%%%%%%%%%%%%%%%%%%%%%%%%%%%%%%%%%%%%%%%%%%%%%%%%%%
%
% Theorem 4.3
\begin{theorem}[Joint Fr\'echet limit]\label{thm:frechet-joint}
%
%%%%%%%%%%%%%%%%%%%%%%%%%%%%%%%%%%%%%%%%%%%%%%%%%%%%%%%%%%%%%%%%%%%%%%%%%%%%%%%%%%%%%%%%%%%%%%%%%%%
Under~\eqref{eq:frechet-tail},
\begin{equation}\label{eq:frechet-joint-limit}
 \left(\frac{R_n}{c_n},\frac{C_n}{c_n}\right)
 \vertk  \bigl(r(\mathcal K),c(\mathcal K)\bigr).
\end{equation}
Since the origin belongs almost surely to the closure of
$\{P_k:k\geq1\}$, it may equivalently be omitted from the definition of
$\mathcal K$.
\end{theorem}
%%%%%%%%%%%%%%%%%%%%%%%%%%%%%%%%%%%%%%%%%%%%%%%%%%%%%%%%%%%%%%%%%%%%%%%%%%%%%%%%%%%%%%%%%%%%%%%%%%%%%%

\smallskip
\begin{proof}
For the compact sample set $\mathcal S_n=\{c_n^{-1}X_1,\ldots,c_n^{-1}X_n\}$,
scaling gives
\[
 r(\mathcal S_n)=\frac{R_n}{c_n},  \qquad  c(\mathcal S_n)=\frac{C_n}{c_n}.
\]
By Lemma~\ref{lem:frechet-compact}, $\mathcal S_n\vertk\mathcal K$ in the Hausdorff metric.  The continuous mapping theorem and Lemma~\ref{lem:hausdorff-stability} therefore give
\eqref{eq:frechet-joint-limit}.
\end{proof}

\medskip

%%%%%%%%%%%%%%%%%%%%%%%%%%%%%%%%%%%%%%%%%%%%%%%%%%%%%%%%%%%%%%%%%%%%%%%%%%%%%%%%%%%%%%%%%%%%%%%%%%%
%
% Corollary 4.4$
%
\begin{corollary}\label{cor:frechet-pareto}
%
%%%%%%%%%%%%%%%%%%%%%%%%%%%%%%%%%%%%%%%%%%%%%%%%%%%%%%%%%%%%%%%%%%%%%%%%%%%%%%%%%%%%%%%%%%%%%%%%%%%%
If, in addition,
\[
 \PP(\rho>s)\sim\lambda s^{-\alpha},
 \qquad s\to\infty,
\]
for some $\lambda>0$, then $c_n\sim(\lambda n)^{1/\alpha}$ and
\[
 \left( (\lambda n)^{-1/\alpha}R_n, (\lambda n)^{-1/\alpha}C_n \right)
 \vertk  \bigl(r(\mathcal K),c(\mathcal K)\bigr).
\]
\end{corollary}

%%%%%%%%%%%%%%%%%%%%%%%%%%%%%%%%%%%%%%%%%%%%%%%%%%%%%%%%%%%%%%%%%%%%%%%%%%%%%%%
%
% Section 5
\section{Joint asymptotics in the Weibull regime}\label{sec:weibull}
%
%%%%%%%%%%%%%%%%%%%%%%%%%%%%%%%%%%%%%%%%%%%%%%%%%%%%%%%%%%%%%%%%%%%%%%%%%%%%%%%

We now turn to the bounded-support case.  Assume throughout this
section that
\begin{equation}\label{eq:weibull-tail}
 0\leq \rho\leq1\quad\text{a.s.},\qquad
 \PP(\rho>1-s)\sim Cs^\gamma,\qquad s\downarrow0,
\end{equation}
where $C>0$ and $\gamma>0$.  Thus, in the notation used in
Henze and Klatt~\cite{hekl}, $\gamma=\beta+1$.  Put
\[
 q_n:=n^{1/\gamma},
 \qquad
 T_{n,j}:=q_n(1-\rho_j),\qquad 1\leq j\leq n.
\]
The marked point process
\[
N_n^{(W)}:=\sum_{j=1}^n \varepsilon_{(T_{n,j},U_j)}
\]
converges, on every bounded height interval $(0,A]$, $A>0$, to a
Poisson point process $N_W$ on $(0,\infty)\times \mathbb{S}^{d-1}$ with intensity
\begin{equation}\label{eq:weibull-intensity}
 C\gamma t^{\gamma-1}\,\D t\,\sigma_{d-1}(\D u).
\end{equation}
Define
\begin{equation}\label{eq:weibull-objectives}
 H_n^{(W)}(z)
 :=
 q_n\left\{
 1-\max_{1\leq j\leq n}\|X_j-q_n^{-1}z\|
 \right\},
 \qquad
 \Psi(z):=\inf_{(t,u)\in N_W}(t+z\cdot u).
\end{equation}
Then
\begin{equation}\label{eq:weibull-value-argmax}
 q_n(1-R_n)=\max_{z\in\mathbb R^d}H_n^{(W)}(z),
 \qquad
 q_nC_n=\argmax_{z\in\mathbb R^d}H_n^{(W)}(z).
\end{equation}

As in the Gumbel regime, the key step is to understand the geometry of the limiting Poisson variational problem.

\medskip

%%%%%%%%%%%%%%%%%%%%%%%%%%%%%%%%%%%%%%%%%%%%%%%%%%%%%%%%%%%%%%%%%%%%%%%%%%%%%%%%%%%%%%%%%%
%
% Lemma 5.1
\begin{lemma}[Geometry of the Weibull limit]\label{lem:weibull-geometry}
%
%%%%%%%%%%%%%%%%%%%%%%%%%%%%%%%%%%%%%%%%%%%%%%%%%%%%%%%%%%%%%%%%%%%%%%%%%%%%%%%%%%%%%%%%%%%%%%%%
Almost surely, $\Psi$ is real-valued, $1$-Lipschitz and concave on
$\mathbb R^d$, satisfies $\Psi(z)\rightarrow-\infty$ as $\|z\|\to\infty$,
and has a unique maximizer $Z_W$.  Writing
\[
 \Phi_W:=\Psi(Z_W),
\]
we have $\Phi_W>0$ almost surely.  Exactly $d+1$ atoms of $N_W$
are active at $Z_W$.  Their directions are affinely independent and
contain the origin in the interior of their convex hull.  In
particular, there is an almost surely positive random variable
$\eta_W$ such that
\begin{equation}\label{eq:weibull-linear-decay}
 \Psi(z)\leq  \Phi_W-\eta_W\|z-Z_W\|,  \qquad z\in\mathbb R^d.
\end{equation}
\end{lemma}
%%%%%%%%%%%%%%%%%%%%%%%%%%%%%%%%%%%%%%%%%%%%%%%%%%%%%%%%%%%%%%%%%%%%%%%%%%%%%%%%%%%%%%%%%%%%%%%%%%%

\smallskip

\begin{proof}
Put $T_*:=\inf_{(t,u)\in N_W} t$.
Then, for $t\ge0$,
$\{T_*>t\} = \big\{N_W\bigl((0,t]\times \mathbb{S}^{d-1}\bigr)=0\big\}$.
Since, by \eqref{eq:weibull-intensity},
\[
\int_{\mathbb{S}^{d-1}}\int_0^t
C\gamma s^{\gamma-1}\,\D s\,\sigma_{d-1}(\D u) = Ct^\gamma,
\]
the void probability formula for Poisson point processes yields
$\PP(T_*>t)=\exp\{-Ct^\gamma\}$, $t\ge0$.
In particular, $0<T_*<\infty$ almost surely.
  On a ball
$\{\|z\|\leq M\}$, an atom with $t>T_*+2M$ cannot affect the
infimum defining $\Psi$.  Indeed,
for $\|z\|\le M$ and any atom $(t,u)\in N_W$,
$t+z\cdot u\ge t-M$.
Since $0<T_*<\infty$ and the process is locally finite, the minimum height is attained almost surely.
Hence, if $(T_*,u_*)$ is an atom with minimal height $T_*$, then
$T_*+z\cdot u_*\le T_*+M$. Thus, whenever $t>T_*+2M$, 
$t+z\cdot u> T_*+M\ge T_*+z\cdot u_*$,
so such an atom cannot attain the infimum defining $\Psi(z)$ on
$\{\|z\|\le M\}$.

Since there are almost surely only finitely
many atoms below any fixed height, $\Psi$ is locally the minimum of
finitely many affine functions.  It is therefore real-valued, concave and
$1$-Lipschitz.

As in the proof of Lemma~\ref{lem:gumbel-geometry}, choose pairwise disjoint open spherical
neighbourhoods of the vertices of a regular simplex centred at the origin
such that, for some $\delta>0$, the convex hull of any choice of one
direction from each neighbourhood contains $\delta\mathbb B^d$.
  Each cap contains almost surely the direction of an atom
of $N_W$, since its intensity integrated over all heights is infinite.
Choose one such atom from each cap, with heights $t_0,\ldots,t_d$.
Then, for every $z$,
\[
 \Psi(z)
 \leq
 \max_{0\leq j\leq d}t_j-\delta\|z\|.
\]
Hence $\Psi(z)\to-\infty$ as $\|z\|\to\infty$ and the maximum is
attained.  Moreover,
\[
 \Phi_W\geq\Psi(0)=T_*>0.
\]

As in the proof of Lemma~\ref{lem:gumbel-geometry}, almost surely every
collection of at most $d$ distinct atom directions is linearly
independent, every collection of $d+1$ distinct directions is
affinely independent, and no $d+2$ affine functions $t+z\cdot u$
associated with distinct atoms take the same value at a common point.
The last assertion follows by conditioning on the directions and the
first $d+1$ heights and using the fact that the remaining height
distribution is non-atomic.

Let $z$ maximize $\Psi$.  In a neighbourhood of $z$, only finitely
many affine functions are relevant. 
If zero did not belong to the convex hull of the directions of the active
functions, the strict separation theorem from convex geometry  would yield a vector
$h\in\mathbb R^d$ and a constant $c>0$ such that
$h\cdot u \ge c$
for every active direction $u$. Hence all active affine functions would
increase under a sufficiently small displacement in the direction $h$.
  Since the inactive functions have a strictly positive gap
above $\Psi(z)$, a sufficiently small displacement in this direction
would increase the minimum, contradicting maximality.  Thus
$ 0\in  \operatorname{conv}\{u:(t,u)\text{ is active at }z\}$.

The general-position properties imply that exactly $d+1$ atoms are
active.  Their directions are affinely independent, and the origin
lies in the interior of their convex hull.

Let $(t_j,u_j)$, $0\leq j\leq d$, be the active atoms.  Since their
convex hull contains a ball $\eta \mathbb{B}^d$ for some $\eta>0$, for every
$h\in\mathbb R^d$, it follows that
$ \min_{0\leq j\leq d}h\cdot u_j \leq-\eta\|h\|$.
Consequently,
\[
 \Psi(z+h)
 \leq
 \min_{0\leq j\leq d}
 \{t_j+(z+h)\cdot u_j\}
 \leq
 \Psi(z)-\eta\|h\|.
\]
This proves uniqueness and~\eqref{eq:weibull-linear-decay}.
\end{proof}

\smallskip

With these geometric properties in hand, the joint Weibull limit follows by the same localization strategy as in the Gumbel case.

\medskip

%%%%%%%%%%%%%%%%%%%%%%%%%%%%%%%%%%%%%%%%%%%%%%%%%%%%%%%%%%%%%%%%%%%%%%%%%%%%%%%%%%%%%%%%%%%%%
%
% Theorem 5.2
%
\begin{theorem}[Joint Weibull limit]\label{thm:weibull-joint}
%%%%%%%%%%%%%%%%%%%%%%%%%%%%%%%%%%%%%%%%%%%%%%%%%%%%%%%%%%%%%%%%%%%%%%%%%%%%%%%%%%%%%%%%%%%%%%%
Under~\eqref{eq:weibull-tail},
\[%begin{equation}\label{eq:weibull-function-convergence}
 H_n^{(W)}\vertk\Psi
 \qquad\text{in }C_{\rm loc}(\mathbb R^d).
\]%end{equation}
Consequently,
\begin{equation}\label{eq:weibull-joint-limit}
 \left(q_n(1-R_n),q_nC_n\right)
 \vertk
 (\Phi_W,Z_W).
\end{equation}
\end{theorem}
%%%%%%%%%%%%%%%%%%%%%%%%%%%%%%%%%%%%%%%%%%%%%%%%%%%%%%%%%%%%%%%%%%%%%%%%%%%%%%%%%%%%%%%%%%%%%%%%%

\smallskip
\begin{proof}

Fix $M>0$ and put $B_M:=\{z\in\mathbb R^d:\|z\|\leq M\}$. For $A>0$, let $N_{n,A}^{(W)}$ and $N_A^{(W)}$ denote the restrictions of $N_n^{(W)}$ and $N_W$, respectively, to $[0,A]\times\mathbb S^{d-1}$.  Put
\[
 \psi_n(t,u;z)  :=  q_n\left\{  1-  \left\|  \left(1-\frac{t}{q_n}\right)u-\frac{z}{q_n}  \right\| \right\}
\]
and define, for $z\in B_M$, 
\[
 \widetilde H_{n,A}^{(W)}(z)  :=  \min\left\{
 A-M,\,  \min_{\substack{1\leq j\leq n\\ T_{n,j}\leq A}}  \psi_n(T_{n,j},U_j;z)  \right\},
\]
as well as
\[
 \widetilde\Psi_A(z)  :=  \min\left\{  A-M,\,  \min_{\substack{(t,u)\in N_W\\ t\leq A}}
 (t+z\cdot u)  \right\}.
\]
The inner minimum is interpreted as $+\infty$ if the corresponding point configuration is empty.  Thus both functions are well defined as elements of $C(B_M)$.

The point-process convergence following from~\eqref{eq:weibull-tail} gives
$ N_{n,A}^{(W)}\vertk N_A^{(W)}$ whenever $N_W$ has no atom at height $A$, which is an event of probability
one.  
The limiting restricted process has only finitely many atoms and is simple.
By Proposition~3.13 of Resnick~\cite{resn}, at every such limiting
configuration the atoms of a vaguely convergent sequence of point measures
can, for all sufficiently large indices, be labelled so that the corresponding
heights and directions converge. Together with the continuous mapping theorem,
this implies the required convergence of the corresponding finite configurations.

Moreover, uniformly for $t\in[0,A]$, $u\in\mathbb S^{d-1}$ and $z\in B_M$,
$\psi_n(t,u;z)=t+z\cdot u+o(1)$,
where the $o(1)$ is deterministic.  Consequently,
$ \widetilde H_{n,A}^{(W)}  \vertk  \widetilde\Psi_A  \qquad\text{in }C(B_M)$.
For $\|z\|\le M$ and $t\in[0,A]$, we have, for all sufficiently large
$n$, $q_n>A$ and hence $1-t/q_n\ge0$. Thus the ordinary and reverse
triangle inequalities give
\[
1-\frac{t}{q_n}-\frac{\|z\|}{q_n}
\le
\left\|\left(1-\frac{t}{q_n}\right)u-\frac{z}{q_n}\right\|
\le
1-\frac{t}{q_n}+\frac{\|z\|}{q_n}.
\]
Consequently,
\begin{equation}\label{eq:weibull-truncation-bound}
 t-M  \leq  \psi_n(t,u;z)  \leq  t+M. 
\end{equation}

Let $T_{n,*}:=\min_{1\leq j\leq n}T_{n,j}$.
On the event $\{T_{n,*}<A-2M\}$, we have $H_n^{(W)}(z)=\widetilde H_{n,A}^{(W)}(z)$, $z\in B_M$.
Indeed, an observation attaining $T_{n,*}$ contributes at most $T_{n,*}+M<A-M$, whereas every observation with height exceeding $A$ contributes more than $A-M$.  Thus observations above the
truncation level cannot determine the minimum, and the artificial level $A-M$ is inactive.  The same argument applies to the limiting process. 
Furthermore, for $A>2M$,
\[
 \PP(T_{n,*}\geq A-2M)  \rightarrow  \exp\{-C(A-2M)^\gamma\},
\]
and the corresponding identity holds for the limiting Poisson process.
Letting $A\to\infty$, the converging-together theorem therefore yields
\[
 H_n^{(W)}|_{B_M}  \vertk  \Psi|_{B_M}  \qquad\text{in }C(B_M).
\]
Since the convergence holds in $C(B_M)$ for every positive integer $M$,
it follows that
\[
H_n^{(W)} \vertk \Psi
\qquad\text{in } C_{\mathrm{loc}}(\mathbb R^d).
\]

Each $H_n^{(W)}$ is real-valued and concave, satisfies $H_n^{(W)}(z)\to - \infty$ as $\|z\| \to \infty$,
and has the unique maximizer $q_nC_n$.  By
Lemma~\ref{lem:weibull-geometry}, the same properties hold for the
limit, with unique maximizer $Z_W$.  Applying
Lemma~\ref{lem:convex-localization} to $-H_n^{(W)}$ and $-\Psi$, and
using~\eqref{eq:weibull-value-argmax}, yields
\eqref{eq:weibull-joint-limit}.
\end{proof}

The same stability mechanism also yields the asymptotic supporting simplex in the Weibull regime.

\medskip

%%%%%%%%%%%%%%%%%%%%%%%%%%%%%%%%%%%%%%%%%%%%%%%%%%%%%%%%%%%%%%%%%%%%%%%%%%%%%%%%%%%%%%%%
%
% Proposition 5.3
\begin{proposition}[Asymptotic supporting simplex]
%
%%%%%%%%%%%%%%%%%%%%%%%%%%%%%%%%%%%%%%%%%%%%%%%%%%%%%%%%%%%%%%%%%%%%%%%%%%%%%%%%%%%%%%%%
\label{prop:weibull-support}
Under~\eqref{eq:weibull-tail},
\[
 \lim_{n\to \infty} \PP(K_n=d+1) = 1.
\]
Moreover, after a suitable ordering, the normalized heights and
directions of the active observations converge jointly to the
$d+1$ active atoms of $N_W$ at $(\Phi_W,Z_W)$.
\end{proposition}
%%%%%%%%%%%%%%%%%%%%%%%%%%%%%%%%%%%%%%%%%%%%%%%%%%%%%%%%%%%%%%%%%%%%%%%%%%%%%%%%%%%%%%%%%%%%%%

\smallskip

\begin{proof}
Fix a ball $B_M:=\{z\in\mathbb R^d:\|z\|\le M\}$
and a bounded height interval $[0,A]$ such that, with probability
arbitrarily close to one, $Z_W$ and $q_nC_n$ belong to $B_M$, the
restriction of $\Psi$ to $B_M$ is determined by the atoms of $N_W$
with heights at most $A$, and, asymptotically, the restriction of
$H_n^{(W)}$ to $B_M$ is determined by the corresponding sample points
with heights at most $A$, as in the proof of
Theorem~\ref{thm:weibull-joint}.

By the Skorokhod representation theorem, as in the proof of
Proposition~\ref{prop:gumbel-support}, we may work on a probability space on which the
restricted finite point measures converge almost surely. Their atoms
may then be matched so that their heights and directions converge.

At the limiting maximizer, Lemma~\ref{lem:weibull-geometry} gives
exactly $d+1$ active atoms; their directions are affinely independent
and contain the origin in the interior of their convex hull.  All other
retained constraints have a strictly positive gap.

For an atom $(T_{n,j},U_{n,j})$ corresponding to the limiting atom
$(t_j,u_j)$, put
\[
 F_{n,j}(z)  :=  q_n\left\{  1-  \left\|  \left(1-\frac{T_{n,j}}{q_n}\right)U_{n,j}
 -\frac{z}{q_n}  \right\|  \right\}.
\]
On a neighbourhood of $Z_W$, the functions $F_{n,j}$ converge in $C^1$
to the affine functions $t_j+z\cdot u_j$.
 Apply
Lemma~\ref{lem:active-stability} to the convex functions $-F_{n,j}$.
The strict gaps exclude all other retained constraints, while
\eqref{eq:weibull-truncation-bound} excludes the discarded ones.
Consequently, for all sufficiently large $n$, precisely the $d+1$
corresponding observations are active at the unique maximizer of
$H_n^{(W)}$.  Letting the probability of the localization complement tend to zero and
returning from the Skorokhod representation to the original probability
space proves the assertion.
\end{proof}

The bounded-support geometry also yields uniform control of the normalized radius deficit and centre. This leads to moment and Wasserstein convergence.

\medskip

%%%%%%%%%%%%%%%%%%%%%%%%%%%%%%%%%%%%%%%%%%%%%%%%%%%%%%%%%%%%%%%%%%%%%%%%%%%%%%%%%%%%%%%%%%%%
%
% Theorem 5.4
%
\begin{theorem}[Uniform tails and convergence of moments]\label{thm:weibull-moments}
%%%%%%%%%%%%%%%%%%%%%%%%%%%%%%%%%%%%%%%%%%%%%%%%%%%%%%%%%%%%%%%%%%%%%%%%%%%%%%%%%%%%%%%%%%%%%%%%
Under~\eqref{eq:weibull-tail}, there exist constants $K_1,K_2>0$,
depending only on $d$ and the distribution of $\rho$, such that
\begin{equation}\label{eq:weibull-tail-bound}
 \PP\left\{ q_n\bigl(1-R_n+\|C_n\|\bigr)>t  \right\}  \leq  K_1 \exp\{-K_2 t^\gamma\},
 \qquad n\geq1,\ t\geq0.
\end{equation}
Consequently, the convergence in
\eqref{eq:weibull-joint-limit} holds in every Wasserstein distance
$W_p$, $1\leq p<\infty$.  More generally, for every continuous
$f:\mathbb R^{d+1}\to\mathbb R$ satisfying
\[
 |f(x)|\leq K(1+\|x\|^p)
\]
for some $K,p<\infty$,
\begin{equation}\label{eq:weibull-moment-convergence}
 \BE f\bigl(q_n(1-R_n),q_nC_n\bigr)  \rightarrow  \BE f(\Phi_W,Z_W).
\end{equation}
Moreover, for all sufficiently small $\lambda>0$,
\begin{equation}\label{eq:weibull-exponential-moments}
 \BE\exp\left\{ \lambda\bigl[q_n(1-R_n+\|C_n\|)\bigr]^\gamma  \right\}
 \rightarrow  \BE\exp\left\{  \lambda(\Phi_W+\|Z_W\|)^\gamma  \right\}  <\infty.
\end{equation}
\end{theorem}
%%%%%%%%%%%%%%%%%%%%%%%%%%%%%%%%%%%%%%%%%%%%%%%%%%%%%%%%%%%%%%%%%%%%%%%%%%%%%%%%%%%%%%%%%%

\smallskip

\begin{proof}
 Put $ \eta:=1/(2\sqrt d)$. For $1\leq j\leq d$, consider spherical caps of Euclidean radius
$\eta$ centred at $\pm e_j$, and let $p_d>0$ denote their common
$\sigma_{d-1}$-measure. For every $c=(c_1,\ldots,c_d)\neq0$, choose an index
$j\in\{1,\ldots,d\}$ such that
$|c_j|=\max_{1\le k\le d}|c_k|$.
Then $|c_j|\ge \|c\|/\sqrt d$.
Take the cap centred at $-\operatorname{sgn}(c_j)e_j$.
 Every $u$ in this cap satisfies
\[
 u\cdot c  \leq  -\frac{\|c\|}{\sqrt d}+\eta\|c\|  =  -\eta\|c\|.
\]
Recall $X_i = \rho_iU_i$. Suppose that every one of the $2d$ caps contains some $U_i$ with
$\rho_i>1-h$. Then, for every $c\in\mathbb R^d$, there exists such an index $i$ for which
\[
 \max_{1\leq i\leq n}\|X_i-c\|  \geq  \|\rho_iU_i-c\|
 \geq  \rho_i-U_i\cdot c  >  1-h+\eta\|c\|.
\]
Taking $c=C_n$ in the preceding inequality and using
$\max_{1\le i\le n}\|X_i-C_n\|=R_n$,
we obtain $1-R_n+\eta\|C_n\|\le h$ and hence
\[
%\begin{equation}\label{eq:weibull-cap-consequence}
 1-R_n+\|C_n\|\leq\frac{h}{\eta},
%\end{equation}
\]
since $1-R_n \ge 0$ and $\eta \le 1$.
Condition~\eqref{eq:weibull-tail} implies that there is $c_0>0$ such
that $\PP(\rho>1-h)\geq c_0h^\gamma$ for $0<h\leq 1$.
Indeed, \eqref{eq:weibull-tail} gives this bound for all sufficiently small $h>0$,
and, after decreasing $c_0$ if necessary, monotonicity extends it
to all $0<h\leq1$.  By independence of $\rho_i$ and $U_i$, for each $i$ the probability that
$U_i$ lies in a prescribed cap and that $\rho_i>1-h$ is at least
$p_dc_0h^\gamma$. Hence the probability that there is no index
$i\in\{1,\ldots,n\}$ with $U_i$ in this cap and $\rho_i>1-h$ is at most
$\exp\{-np_dc_0h^\gamma\}$. 
A union bound over the $2d$ caps therefore yields
\[
 \PP\left\{ 1-R_n+\|C_n\|>\frac{h}{\eta} \right\} \leq 2d\exp\{-np_dc_0h^\gamma\}.
\]
For $t\leq q_n/\eta$, take $h=\eta t/q_n$.  Since $q_n^\gamma=n$,
\[
 \PP\left\{ q_n(1-R_n+\|C_n\|)>t \right\} \leq 2d\exp\{-p_dc_0\eta^\gamma t^\gamma\}.
\]
If $t>q_n/\eta$, the probability on the left is zero.  Indeed, the
centre of a smallest enclosing ball belongs to the convex hull of the
sample, so $\|C_n\|\leq1$, and $R_n\leq1$, whence
$ 1-R_n+\|C_n\|\leq2\leq 1/\eta$.
This proves~\eqref{eq:weibull-tail-bound}.

The bound implies uniform integrability of every positive power of
$\bigl\|\bigl(q_n(1-R_n),q_nC_n\bigr)\bigr\|$.
Together with the weak convergence in~\eqref{eq:weibull-joint-limit},
this yields convergence of the corresponding $p$th moments for every
$p\ge1$, and hence convergence in the Wasserstein distance $W_p$.
More generally, the same uniform integrability argument gives
\eqref{eq:weibull-moment-convergence} for every continuous function
$f:\mathbb R^{d+1}\to\mathbb R$ of polynomial growth.

Finally, \eqref{eq:weibull-tail-bound} implies uniform integrability of
$ \exp\left\{ \lambda[q_n(1-R_n+\|C_n\|)]^\gamma \right\}$ 
whenever $\lambda>0$ is sufficiently small.  The continuous mapping
theorem and uniform integrability then yield
\eqref{eq:weibull-exponential-moments}.
\end{proof}

We next give an explicit description of the joint law of the limiting
value and optimizer.  To this end, put $m=d+1$.  For
$u=(u_0,\ldots,u_d)\in(\mathbb S^{d-1})^m$, define
\[
 \Delta(u)  :=  \left| \det 
 \begin{pmatrix} 1&u_0^\top\\
 \vdots&\vdots\\
 1&u_d^\top
 \end{pmatrix}
 \right|,  \qquad  B(u)  :=  \mathbf 1\{0\in\operatorname{int}  \operatorname{conv}(u_0,\ldots,u_d)\}.
\]
For $\gamma>0$, let $g_\gamma(v) := v^{\gamma-1}$ if $v>0$ and $g_\gamma(v) := 0$, otherwise.
 Moreover, putting $t_+:=\max\{t,0\}$ for $t \in \mathbb{R}$, define
 \[
%\begin{equation}\label{eq:Agamma}
 A_\gamma(y)
 :=  \int_{\mathbb S^{d-1}}  (1-y\cdot u)_+^\gamma\,\sigma_{d-1}(\D u)
%\end{equation}
\]
and
\[
%\begin{equation}\label{eq:Jgamma}
 J_\gamma(y)
 :=  \int_{(\mathbb S^{d-1})^m}  \Delta(u)B(u)  \prod_{i=0}^d  g_\gamma(1-y\cdot u_i)\,
 \sigma_{d-1}^{\otimes m}(\D u).
%\end{equation}
\]
When $\gamma<1$, the value assigned to $g_\gamma(0)$ is immaterial;
all density identities below are understood to hold Lebesgue-almost everywhere.

Beyond weak and moment convergence, the Poisson representation allows the limiting value--optimizer pair to be described explicitly.

\medskip

%%%%%%%%%%%%%%%%%%%%%%%%%%%%%%%%%%%%%%%%%%%%%%%%%%%%%%%%%%%%%%%%%%%%%%%%%%%%%%%%%%%%%%%%%%%%%%%%%%
%
% Theorem 5.5
\begin{theorem}[Joint density and gamma factorisation]\label{thm:weibull-density}
%
%%%%%%%%%%%%%%%%%%%%%%%%%%%%%%%%%%%%%%%%%%%%%%%%%%%%%%%%%%%%%%%%%%%%%%%%%%%%%%%%%%%%%%%%%%%%%%%
The pair $(\Phi_W,Z_W)$ has density
\begin{align}
 q_W(\phi,z)
 ={}&  \frac{(C\gamma)^{d+1}}{(d+1)!}
 \exp\left\{ -C\int_{\mathbb S^{d-1}}        (\phi-z\cdot u)_+^\gamma\,        \sigma_{d-1}(\D u)
 \right\}
 \nonumber\\
 &\quad\times  \int_{(\mathbb S^{d-1})^{d+1}}  \Delta(u)B(u)  \prod_{i=0}^d
 g_\gamma(\phi-z\cdot u_i)\,  \sigma_{d-1}^{\otimes(d+1)}(\D u),  \label{eq:weibull-joint-density}
\end{align}
for $\phi>0$ and $z\in\mathbb R^d$. Moreover, putting
\begin{equation}\label{eq:weibull-QY}
 Y_W:=\frac{Z_W}{\Phi_W},  \qquad  Q_W:=C\Phi_W^\gamma A_\gamma(Y_W),
\end{equation}
then $Q_W$ and $Y_W$ are independent, with
\begin{equation}\label{eq:weibull-gamma}
 Q_W\sim\operatorname{Gamma}(d+1,1),
\end{equation}
while $Y_W$ has density
\begin{equation}\label{eq:YW-density}
 f_{Y_W}(y)  =  \frac{\gamma^d}{d+1}  \frac{J_\gamma(y)}{A_\gamma(y)^{d+1}},  \qquad y\in\mathbb R^d.
\end{equation}
Consequently,
\begin{equation}\label{eq:weibull-representation}
 (\Phi_W,Z_W)  \edist  \left[   \frac{Q_W}{C A_\gamma(Y_W)}  \right]^{1/\gamma}  (1,Y_W),
\end{equation}
where $Q_W$ and $Y_W$ on the right-hand side are independent.
\end{theorem}

\smallskip

\begin{proof}
By Lemma~\ref{lem:weibull-geometry}, almost surely the optimum is unique
and is determined by a unique unordered set of $d+1$ active atoms.
Hence the corresponding ordered $(d+1)$-tuples are counted exactly
$(d+1)!$ times. Apply the multivariate Mecke formula
(see, e.g., Last and Penrose~\cite[Theorem~4.4]{LP17})
to ordered $(d+1)$-tuples of atoms and divide by $(d+1)!$.
For fixed directions $u_0,\ldots,u_d$ such that
\[
0\in\operatorname{int}\operatorname{conv}\{u_0,\ldots,u_d\},
\]
an active tuple corresponding to value $\phi$ and optimizer $z$ has heights
$ t_i=\phi-z\cdot u_i$ for $0\leq i\leq d$.
The absolute Jacobian of the transformation $ (\phi,z)\mapsto(t_0,\ldots,t_d)$ 
is $\Delta(u)$.  For this tuple to determine the optimum, every
remaining Poisson point must satisfy
$t+z\cdot u\geq\phi$.

Conversely, if this void condition holds and $B(u)=1$, then $z$ is
the unique global maximizer. Indeed, since
$0\in\operatorname{int}\operatorname{conv}\{u_0,\ldots,u_d\}$,
there exists $\eta>0$ such that, for every $h\in\mathbb R^d$,
\[
 \min_{0\le i\le d} h\cdot u_i \le -\eta\|h\|.
\]
Hence $ \Psi(z+h) \le  \min_{0\le i\le d}\{t_i+(z+h)\cdot u_i\}  \le  \phi-\eta\|h\|$,
whereas the void condition gives $\Psi(z)=\phi$.

The intensity of the complementary set is
\[
 \int_{\mathbb S^{d-1}}    \int_0^{(\phi-z\cdot u)_+}    C\gamma t^{\gamma-1}\,\D t\,
   \sigma_{d-1}(\D u)   =  C\int_{\mathbb S^{d-1}}    (\phi-z\cdot u)_+^\gamma\,    \sigma_{d-1}(\D u).
\]
Thus the corresponding void probability is the exponential factor in
\eqref{eq:weibull-joint-density}.  Multiplying it by the $d+1$
height-intensity densities, the Jacobian and the factor
$1/(d+1)!$ gives~\eqref{eq:weibull-joint-density}.

For the factorisation, set $z=\phi y$.  Since
$ (\phi-\phi y\cdot u)_+^\gamma  =  \phi^\gamma(1-y\cdot u)_+^\gamma$
and
\[
 g_\gamma(\phi-\phi y\cdot u)  =  \phi^{\gamma-1}g_\gamma(1-y\cdot u),
\]
while the transformation $(\phi,y)\mapsto(\phi,z)$ has Jacobian
$\phi^d$, the joint density of $(\Phi_W,Y_W)$ is
\begin{equation}\label{eq:phi-Y-density}
 \frac{(C\gamma)^{d+1}}{(d+1)!}  J_\gamma(y)  \phi^{(d+1)\gamma-1}  \exp\{-C\phi^\gamma A_\gamma(y)\}.
\end{equation}
Now put $q=C\phi^\gamma A_\gamma(y)$. A direct change of variables transforms
\eqref{eq:phi-Y-density} into
\[
 \frac{q^d\E^{-q}}{d!}
 \,  \frac{\gamma^d}{d+1}  \frac{J_\gamma(y)}{A_\gamma(y)^{d+1}},  \qquad q>0.
\]
This factorisation proves
\eqref{eq:weibull-gamma}--\eqref{eq:YW-density}, and \eqref{eq:weibull-representation} follows from
\eqref{eq:weibull-QY}.
\end{proof}

\smallskip

In the uniform-ball case, this Poisson variational problem admits a particularly transparent geometric interpretation.

\medskip

%%%%%%%%%%%%%%%%%%%%%%%%%%%%%%%%%%%%%%%%%%%%%%%%%%%%%%%%%%%%%%%%%%%%%%%%%%%%%%%%%%%%%%%%%%%%%%%%%%%
%
% Corollary 5.6
\begin{corollary}[Uniform ball and the Poisson zero cell]\label{cor:uniform-ball-zero-cell}
%
%%%%%%%%%%%%%%%%%%%%%%%%%%%%%%%%%%%%%%%%%%%%%%%%%%%%%%%%%%%%%%%%%%%%%%%%%%%%%%%%%%%%%%%%%%%%%%%%%%
Suppose that $X$ is uniformly distributed on the unit ball $\mathbb{B}^d$.
Then $\gamma=1$, $C=d$ and $q_n=n$.  Let
\[
 Z:=  \bigcap_{(t,u)\in N_W}  \{x\in\mathbb R^d:x\cdot u\leq t\}.
\]
If $r_{\rm in}(Z)$ and $c_{\rm in}(Z)$ denote, respectively, the
inradius and the almost surely unique incenter of $Z$,
then almost surely \begin{equation}\label{eq:zero-cell-identification}
 (\Phi_W,Z_W)  =  \bigl(r_{\rm in}(Z),-c_{\rm in}(Z)\bigr).
\end{equation}
Consequently, $ \bigl(n(1-R_n),nC_n\bigr)  \vertk  \bigl(r_{\rm in}(Z),-c_{\rm in}(Z)\bigr)$.
\end{corollary}
%%%%%%%%%%%%%%%%%%%%%%%%%%%%%%%%%%%%%%%%%%%%%%%%%%%%%%%%%%%%%%%%%%%%%%%%%%%%%%%%%%%%%%%%%%%%%%%%

\smallskip

\begin{proof}
For the uniform distribution on $\mathbb{B}^d$, we have
$\PP(\rho>1-s)  =  1-(1-s)^d  \sim ds$,
so that $\gamma=1$ and $C=d$.  Moreover,
$ d\,\D t\,\sigma_{d-1}(\D u)  =  \kappa_d^{-1}\,\D t\,\mathcal H^{d-1}(\D u)$,
which is the isotropic hyperplane intensity appearing in the
translation--scaling specialization of
Kabluchko, Marynych and Molchanov~\cite{KMM}.

For $q\geq0$, $\Psi(z)\geq q$ is equivalent to $ t+z\cdot u\geq q$ for every $(t,u)\in N_W$,
or, equivalently, $(-z)\cdot u+q\leq t$ for every $(t,u)\in N_W$. The latter condition is exactly
$-z+q\mathbb{B}^d\subset Z$. Hence
\[
 \{z:\Psi(z)\geq q\}  =  -\{x:x+q\mathbb{B}^d\subset Z\}.
\]
Maximizing over $q$ and using the almost sure uniqueness of the
maximizer from Lemma~\ref{lem:weibull-geometry} proves
\eqref{eq:zero-cell-identification}.
\end{proof}

\medskip

%%%%%%%%%%%%%%%%%%%%%%%%%%%%%%%%%%%%%%%%%%%%%%%%%%%%%%%%%%%%%%%%%%%%%%%%%%%%%%%%%%%%%%%%%%%%%
%
% Remark 5.7
\begin{remark}\label{rem:KMM-comparison}
%
%%%%%%%%%%%%%%%%%%%%%%%%%%%%%%%%%%%%%%%%%%%%%%%%%%%%%%%%%%%%%%%%%%%%%%%%%%%%%%%%%%%%%%%%%
Corollary~\ref{cor:uniform-ball-zero-cell} makes precise the connection
with the random-set limit of Kabluchko, Marynych and
Molchanov~\cite{KMM}.  Their limit theorem concerns the entire
rescaled feasible set of affine transformations.  Restriction to
translations and scalar dilations yields the half-space constraints
above.  The joint radius--centre limit additionally uses localization
and the almost sure uniqueness of the optimizing pair.  Theorem
\ref{thm:weibull-joint} extends this uniform-ball situation to the
general boundary behaviour~\eqref{eq:weibull-tail}.
\end{remark}

\medskip

\section*{Acknowledgement}
I am especially grateful to Ilya Molchanov for drawing my attention to the work of Kabluchko, Marynych and Molchanov, and, more importantly, for sharing his geometric insight into the problem. His perspective substantially influenced the way I came to view the problem and was instrumental in shaping the present approach.

\section*{A note on the use of AI}
During the preparation of this paper, the author used ChatGPT (OpenAI, GPT-5.6 Sol) as an auxiliary tool for language editing, presentation, and the critical review of mathematical arguments. All suggestions were independently checked by the author, who assumes full responsibility for the mathematical content and final formulation.

\end{document}